
%


\def\juerg #1{{\color{red}#1}}
\let\juerg\relax

\def\revis#1{{\color{red}#1}}
\def\revis#1{#1}


\def\LaTeX{%
  \let\Begin\begin
  \let\End\end
  \def\Bcenter{\Begin{center}}
  \def\Ecenter{\End{center}}
  \let\Label\label
  \let\salta\relax
  \let\finqui\relax
  \let\futuro\relax}

\def\UK{\def\our{our}\let\sz s}
\def\USA{\def\our{or}\let\sz z}



\LaTeX

\USA


\salta

\documentclass[twoside,12pt]{article}
\setlength{\textheight}{24cm}
\setlength{\textwidth}{16cm}
\setlength{\oddsidemargin}{2mm}
\setlength{\evensidemargin}{2mm}
\setlength{\topmargin}{-15mm}
\parskip2mm


%
%
\usepackage{cite}

\usepackage{color}
\usepackage{amsmath}
\usepackage{amsthm}
\usepackage{amssymb}

\usepackage{amsfonts}
\usepackage{mathrsfs}

\usepackage{hyperref}
\usepackage[mathcal]{euscript}

\usepackage[ulem=normalem,draft]{changes}




%
\newtheorem{theorem}{Theorem}[section]

\newtheorem{corollary}[theorem]{Corollary}

\newtheorem{lemma}[theorem]{Lemma}

\newenvironment{proofteo1}{\noindent{{\it Proof of Theorem~\ref{Teo1}.}}}{\hfill$\square$}

\finqui

\def\Beq{\Begin{equation}}
\def\Eeq{\End{equation}}
\def\Bsist{\Begin{eqnarray}}
\def\Esist{\End{eqnarray}}

\def\Bthm{\Begin{theorem}}
\def\Ethm{\End{theorem}}
\def\Blem{\Begin{lemma}}
\def\Elem{\End{lemma}}

\def\Bcor{\Begin{corollary}}
\def\Ecor{\End{corollary}}
\def\Brem{\Begin{remark}\rm}
\def\Erem{\End{remark}}

\def\Bdim{\Begin{proof}}
\def\Edim{\End{proof}}
\let\non\nonumber




\def\step #1 \par{\medskip\noindent{\bf #1.}\quad}


\def\Lip{Lip\-schitz}

\def\lhs{left-hand side}
\def\rhs{right-hand side}



\def\multibold #1{\def\arg{#1}%
  \ifx\arg\pto \let\next\relax
  \else
  \def\next{\expandafter
    \def\csname #1#1#1\endcsname{{\bf #1}}%
    \multibold}%
  \fi \next}

\def\pto{.}

\def\multical #1{\def\arg{#1}%
  \ifx\arg\pto \let\next\relax
  \else
  \def\next{\expandafter
    \def\csname cal#1\endcsname{{\cal #1}}%
    \multical}%
  \fi \next}


\def\multimathop #1 {\def\arg{#1}%
  \ifx\arg\pto \let\next\relax
  \else
  \def\next{\expandafter
    \def\csname #1\endcsname{\mathop{\rm #1}\nolimits}%
    \multimathop}%
  \fi \next}

\multibold
qwertyuiopasdfghjklzxcvbnmQWERTYUIOPASDFGHJKLZXCVBNM.

\multical
QWERTYUIOPASDFGHJKLZXCVBNM.

\multimathop
dist div dom meas sign supp .


\def\Accorpa #1#2 #3 {\gdef #1{\eqref{#2}--\eqref{#3}}%
  \wlog{}\wlog{\string #1 -> #2 - #3}\wlog{}}


\def\<#1>{\mathopen\langle #1\mathclose\rangle}
\def\norma #1{\mathopen \| #1\mathclose \|}
\def\normaH #1{\mathopen \| #1\mathclose \|_H}


\def\separa{\noalign{\allowbreak}}

\def\infess{\mathop{\rm ess\,inf}}
\def\supess{\mathop{\rm ess\,sup}}

\def\iO{\int_\Omega}
\def\iQt{\iint_{Q_t}}
\def\iQ{\iint_Q}

\def\dt{\partial_t}
\def\dn{\partial_{\bf n}}

\def\checkmmode #1{\relax\ifmmode\hbox{#1}\else{#1}\fi}


\def\erre{{\mathbb{R}}}
\def\enne{{\mathbb{N}}}




\def\genspazio #1#2#3#4#5{#1^{#2}(#5,#4;#3)}
\def\spazio #1#2#3{\genspazio {#1}{#2}{#3}T0}

\def\L {\spazio L}
\def\H {\spazio H}
\def\W {\spazio W}

\def\C #1#2{C^{#1}([0,T];#2)}


\def\Lx #1{L^{#1}(\Omega)}
\def\Hx #1{H^{#1}(\Omega)}

\def\Luno{\Lx 1}
\def\Ldue{\Lx 2}

\def\Huno{\Hx 1}
\def\Hdue{\Hx 2}

\def\LiQ{L^\infty(Q)}


\def\LQ #1{L^{#1}(Q)}


\let\theta\vartheta
\let\eps\varepsilon
\let\phi\varphi

\let\TeXchi\chi                         
\newbox\chibox
\setbox0 \hbox{\mathsurround0pt $\TeXchi$}
\setbox\chibox \hbox{\raise\dp0 \box 0 }
\def\chi{\copy\chibox}



\def\Betaeps{f_{1,\eps}}
\def\betaeps{f'_{1,\eps}}
\def\phieps{\phi_\eps}
\def\mueps{\mu_\eps}
\def\phiz{\phi_0}


\def\CX{{\cal X}}

\def\CS{{\cal S}}

\def\CY{{\cal Y}}
\def\CU{{\cal U}}

\def\Uad{{\cal U}_{\rm ad}}
\def\CZ{{\cal Z}}
\def\CN{{\cal N}}

\def\us{u^*}
\def\phis{\phi^*}
\def\mus{\mu^*}
\def\ws{w^*}
\def\qs{q^*}
\def\ps{p^*}
\def\rs{r^*}

\def\by{{\bf y}}
\def\bys{{\bf y^*}}

\def\xih{\xi^h}
\def\xik{\xi^k}

\def\VD{V^*}


\let\tau \tau
\normalfont

\Begin{document}


\title{{\bf Optimality conditions for sparse\\ optimal control of viscous\\ Cahn--Hilliard systems with
\\logarithmic potential}}
\author{}
\date{}
\maketitle

\Bcenter
\vskip-1.9cm
{\large\bf Pierluigi Colli$^{(1)}$}\\
{\normalsize e-mail: {\tt pierluigi.colli@unipv.it}}\\[.4cm]
{\large\bf J\"urgen Sprekels$^{(2)}$}\\
{\normalsize e-mail: {\tt juergen.sprekels@wias-berlin.de}}\\[.4cm]
{\large\bf Fredi Tr\"oltzsch$^{(3)}$}\\
{\normalsize e-mail: {\tt  troeltzsch@math.tu-berlin.de}}\\[.6cm]

$^{(1)}$
{\small Dipartimento di Matematica ``F. Casorati'', Universit\`a di Pavia}\\
{\small and Research Associate at the IMATI -- C.N.R. Pavia}\\ 
{\small via Ferrata 5, 27100 Pavia, Italy}\\[.2cm]
$^{(2)}$
{\small Weierstrass Institute for Applied Analysis and Stochastics}\\
{\small Mohrenstra\ss e 39, 10117 Berlin, Germany}\\[.2cm]
$^{(3)}$
{\small Institut f\"ur Mathematik der Technischen Universit\"at Berlin}\\
{\small Stra\ss e des 17.~Juni 136, 10623 Berlin, Germany}\\[.8cm]
\Ecenter

{
\Begin{abstract}\noindent
In this paper we study the optimal control of a parabolic initial-boundary value problem
 of viscous Cahn--Hilliard type with zero Neumann boundary
conditions. Phase field systems of this type govern the evolution of
diffusive phase transition processes with conserved order parameter. 
It is assumed that the nonlinear \revis{functions} driving the physical processes within the spatial domain are 
double-well potentials of logarithmic type whose derivatives become singular at the boundary of their respective domains of definition. 
For such systems, optimal control problems have been studied in the past.  We focus here on the situation when the cost functional 
of the optimal control problem contains a nondifferentiable term like the $L^1$-norm, which leads to sparsity of optimal controls. 
For such cases, we establish first-order necessary and second-order sufficient optimality conditions for locally optimal controls.
In the approach to second-order sufficient conditions, the main novelty of this paper, we adapt a technique introduced 
by E.~Casas, C.~Ryll and F.~Tr\"{o}ltzsch  in the paper 
 [{\em SIAM J. Control Optim.} {\bf  53} (2015), 2168--2202].
In this paper, we show that this method can also be successfully applied 
to systems of viscous Cahn--Hilliard type with logarithmic nonlinearity. Since the Cahn--Hilliard system corresponds to a 
fourth-order partial differential equation in contrast to the second-order systems investigated before, additional technical difficulties 
have to be overcome.      
\\[2mm]
\normalfont
\\[2mm]

\End{abstract}
}
\salta

\pagestyle{myheadings}
\newcommand\testopari{\sc Colli \ --- \ Sprekels \ --- \ Tr\"oltzsch}
\newcommand\testodispari{\sc Optimality conditions with sparsity for a Cahn--Hilliard system}
\markboth{\testopari}{\testodispari} 


\finqui


\section{Introduction}
\label{Intro}
\setcounter{equation}{0}

Let $\Omega\subset \erre^3$ denote some bounded and connected open set with smooth boundary $ \Gamma=\partial\Omega$ (a 
compact hypersurface of class $C^2$) and unit outward 
normal ${\bf n}$.
Moreover, let $T>0$ denote some final time, and set
\begin{align*}
&Q_t:=\Omega \times (0,t),\quad \Sigma_t:=\Gamma\times (0,t), \quad \mbox{for }\,t\in(0,T],\quad\mbox{and}\quad
Q:=Q_T,\quad \Sigma:=\Sigma_T.
\end{align*}
We then study the following optimal control problem:

\vspace*{2mm}\noindent
{\bf (CP)} \,\,Minimize the cost functional
\begin{align}
\label{cost}
{\cal J}(\phi, u)\,&:=\,\frac{b_1}2\iQ|\phi -\phi_Q|^2\,+\,\frac{b_2}2 \iO|\phi(T)-\phi_{\Omega}|^2\,+\,\frac{b_3}2\iQ|u|^2
\,+\,\kappa\,\juerg{G}(u)\,,\nonumber\\
&=: J(\phi,u)\,+\,\kappa\,\juerg{G}(u)
\end{align}
subject to the initial-boundary value system 
\begin{align}
\label{ss1}
&\dt \phi -\Delta \mu = 0 &&\mbox{a.e. in }\,Q,\\
\label{ss2}
&\tau \dt \phi -\Delta \phi + f'(\phi) = \mu + w  &&\mbox{a.e. in }\,Q,\\
\label{ss3}
&\gamma \dt w + w = u   &&\mbox{a.e. in }\,Q,\\
\label{ss4}
&\dn \mu = \dn\phi = 0 &&\mbox{a.e. on }\,\Sigma,\\
\label{ss5}
&\phi(0)=\phi_0,\,\quad w(0)= w_0 &&\mbox{a.e. in }\,\Omega,
\end{align}
\Accorpa\State ss1 ss5
and to the control constraint
\begin{equation}
\label{defUad}
\Uad=\{u\in\CU: \ \underline u(x,t)\le u(x,t)\le \overline u(x,t) \,\mbox{ for a.a. $(x,t)$\, in }\,Q\}.
\end{equation}
Here, the given bounds $\,\underline u,\overline u\in \CU\,$ satisfy $\,\underline u\le\overline u$ almost everywhere in $Q$, and
 the control space is given by
\begin{equation}
\label{defU}
\CU=\LiQ.
\end{equation}
Moreover, the targets $\phi_Q, \, \phi_{\Omega}$ are given functions, \juerg{$b_1\ge 0$, $b_2\ge 0$, $b_3>0$ are constants,}
and
$\kappa>0$ is a constant which represents the sparsity parameter. The sparsity-enhancing functional $\,\juerg{G}:L^2(Q)\to\erre\,$ is
\juerg{nonnegative,} continuous and convex. Typically, \juerg{$G$} has a nondifferentiable form like, e.g.,
\begin{equation}
\label{defg}
\juerg{G}(u)=\|u\|_{L^1(Q)}=\iint_Q|u|\,. 
\end{equation}

The state equations \eqref{ss1}--\eqref{ss2} constitute a viscous Cahn--Hilliard system (introduced in \cite{CH}), in which a number of 
physical constants have been normalized to unity and whose state variables 
$\,\varphi\,$ and $\,\mu\,$  are monitored through 
the input variable $\,w$, which is in turn determined by the action of the control $\,u\,$ via the linear control equation~\eqref{ss3}. 
Equation~\eqref{ss3} models how the ``forcing'' $w$ is generated by the external control $u$. We remark that \eqref{ss4}
could be replaced by much more general differential equations modeling the relation between an $L^2$-control $u$ and a smooth forcing $w$: 
one can see, e.g, the system studied in \cite{CGRS}, in which the linear equation~\eqref{ss3} is replaced by a reaction-diffusion equation where the unknown $w$ represents a nutrient concentration, in a model for tumor growth.

In the system \eqref{ss1}--\eqref{ss5}, $\,\varphi\,$ plays the role of an \emph{order parameter} that attains its values in the
interval $[-1,+1]$, while $\,\mu\,$ is the associated \emph{chemical potential}.  
Moreover,  $\tau>0 $ is a viscosity coefficient,  $\gamma $  is a given (uniformly) positive function defined on $\Omega$, and
$ \phi_0$ and $ w_0$ are given initial data. The nonlinearity  $\,f\,$ represents a double-well potential whose
derivative defines the local part of the thermodynamic force driving the evolution of the system. 
In this paper, we consider potentials having the typical behavior of the physically particularly relevant 
logarithmic potential given by
\begin{align}
\label{flog}
f_{\rm log}(r)&=\left\{\begin{array}{ll}
c_1\bigl((1+r)\,\ln(1+r)+(1-r)\,\ln(1-r)\bigr)-c_2r^2 &\quad\mbox{if }\,r\in(-1,1)\\
2c_1\ln(2)-c_2 &\quad\mbox{if }\,r\in\{-1,1\}\\
+\infty&\quad\mbox{if }\,r\not\in [-1,1].
\end{array}
\right. 
\end{align}
In this connection, $c_1,c_2$ are nonnegative and such that $f_{\rm log}$ is nonconvex. Notice that 
for $f=f_{\rm log}$ the term $f'(\phi)$
occurring in \eqref{ss2} becomes singular as $\phi\searrow -1$ and $\phi\nearrow1$, which forces the order parameter
$\phi$ to attain its values in the physically meaningful range~$(-1,1)$.  

Starting with the seminal paper \cite{EZ}, there exists an abundant literature on the well-posedness and asymptotic behavior of viscous and nonviscous Cahn--Hilliard
systems with zero Neumann and with dynamic boundary conditions that cannot be cited here in its entirety. A nice 
collection of papers on this topic 
up to the year 2015 can be found in \cite{Heida}. In spite of this large amount of related literature, we have chosen
to provide a detailed well-posedness analysis of the state system \State, both for the readers' convenience and the fact that 
the system \State\ was apparently not studied before in this particular form in which the control contributes to the chemical potential
through the quantity $w$. Notice that the typical regularity to be expected for an $L^2-$control $u$ is $w\in H^1(0,T;\Ldue)$, which
in the three-dimensional case with logarithmic potential  is typically needed to derive a {\em separation property} 
from \eqref{ss2} for the state variable $\varphi$.
   
There also exist contributions to the optimal control of Cahn--Hilliard type systems in various contexts. 
Without claiming to be exhaustive and complete, we mention now some related papers. First, let us
refer to~\cite{Duan,HintWeg,Z,ZW} and, in the framework of diffusive models of tumor growth, to 
\cite{CGRS, CSig,CSS1,CSS2,EK1,EK2,GLR}. Problems with dynamical boundary conditions have been studied in  
\cite{CFGS1,CFGS2,CGSANA,CGSAMO,CGSAnnali,CGSSIAM,CGSconv,CSig,GS}, and convective Cahn--Hilliard systems have been the subject of
\cite{CGSAnnali,CGSSIAM,GS,RoSp,ZL1,ZL2}. In addition, quite a number of works have been dedicated to the study of cases in which 
the Cahn--Hilliard system is coupled to other systems; in this connection, we quote Cahn--Hilliard--Navier--Stokes models
(see \cite{FGS,HKW,HW2,HW3,Medjo}) and the Cahn--Hilliard--Oono (see \cite{CGRS2,GiRoSi}),  Cahn--Hilliard--Darcy 
(see\cite{ACGW,SpWu}), Cahn--Hilliard--Brinkman (see~\cite{EK1}) and Cahn--Hilliard with curvature effects (see~\cite{CGSS6}) systems.

None of the papers cited above is concerned with the aspect of {\em sparsity}, i.e., the 
possibility that any locally optimal control may vanish in subregions of positive measure of the space-time cylinder $Q$
 that are controlled by the
sparsity parameter $\kappa$. The geometry of these subregions depends on the particular choice of the convex functional $\juerg{G}$,
which can differ in different situations. The sparsity properties can be deduced from the variational inequality occurring in the first-order
necessary optimality conditions and the particular
form of the subdifferential  $\partial \juerg{G}$. 
In this paper, we focus on sparsity, where, in the following, we restrict ourselves to the case of {\em full sparsity}
which is connected to the $L^1(Q)$-norm functional \juerg{$G$} introduced in~\eqref{defg}. Other types of sparsity such as 
{\em directional sparsity with respect to time} 
and {\em directional sparsity with respect to space} (see, e.g.,~\cite{SpTr1}) are not treated in this paper.

Sparsity in the optimal control theory for partial differential equations has become an actively investigated aspect.
The use of sparsity-enhancing functionals goes back to inverse problems and image processing. It was the seminal paper  
\cite{stadler2009}  on elliptic control problems that initiated the discussion of sparsity in the optimal control theory 
of partial differential equations. Soon after  \cite{stadler2009}, many results on sparse optimal controls for PDEs 
were published. We mention only very few of them with closer relation to our paper, in particular 
\cite{casas_herzog_wachsmuth2017,herzog_obermeier_wachsmuth2015,herzog_stadler_wachsmuth2012}, on directional sparsity,
 and \cite{casas_troeltzsch2012} on a general theorem for second-order conditions. 
Moreover, we refer to some new trends in the investigation of sparsity, namely, infinite horizon sparse 
optimal control (see, e.g., \cite{Kalise_Kunisch_Rao2017,Kalise_Kunisch_Rao2020}) and fractional order optimal 
control (cf.~\cite{Otarola2020,Otarola_Salgado2018}).

The abovementioned papers concentrated on the first-order optimality conditions for sparse optimal controls of single elliptic and parabolic equations. 
In  \cite{casas_ryll_troeltzsch2013,casas_ryll_troeltzsch2015}, first- and second-order optimality conditions 
have been discussed in the context of sparsity for the (semilinear) system of  FitzHugh--Nagumo equations. 
More recently, sparsity of optimal controls for reaction-diffusion systems of Cahn--Hilliard type have been addressed in 
\cite{CSS4,Garcke_etal2021,SpTr1}.  Moreover, we refer to the measure control of the Navier--Stokes system 
studied in \cite{Casas_Kunisch2021}.
However, to the best knowledge of the authors, second-order sufficient optimality for sparse controls 
for the Cahn--Hilliard and viscous Cahn--Hilliard equations have never been studied before. 

Second-order sufficient optimality 
conditions are usually based on a condition of coercivity that is required to hold for the smooth part  $\,J\,$ 
of $\,{\cal J}\,$ in a certain {\em critical cone}. The nonsmooth part $\,\juerg{G}\,$ contributes to sufficiency by its convexity. 
For the strength of sufficient conditions it is crucial that the critical cone be as small as possible. In their paper
\cite{casas_ryll_troeltzsch2015}, Casas--Ryll--Tr\"oltzsch devised a technique by means of which a very advantageous (i.e., small)
critical cone can be chosen. This method was originally introduced  for a class of semilinear second-order parabolic problems
with smooth nonlinearities. In the recent papers \cite{SpTr2,SpTr3} two of the present authors have demonstrated that 
it can be adapted correspondingly to the sparse optimal control of Allen--Cahn systems
with dynamic boundary conditions and to a large class of systems modeling tumor growth, where in both papers the
case of singular logarithmic nonlinearities of the form \eqref{flog} was admitted.

It is the main aim and novelty of this work to show
that also systems having a Cahn--Hilliard structure can be treated accordingly (at least in the viscous case $\,\tau>0$). This
is by no means obvious, since, in contrast to the second-order systems investigated in \cite{SpTr2,SpTr3}, the Cahn--Hilliard structure
studied here leads to a fourth-order PDE for the order parameter
 $\varphi$ (which readily follows from insertion for $\mu$ from \eqref{ss2} in \eqref{ss1}). As a consequence, a number of 
additional technical difficulties have to be overcome, both in the proof of the Fr\'echet differentiability of the control-to-state 
operator and in the analysis of the properties of the adjoint variables. Some of these technical difficulties are also 
due to the singular behavior of the derivative \,$f'(\varphi)$\, of the logarithmic nonlinearity appearing in \eqref{ss2}. The nonviscous case
$\tau=0$ with logarithmic nonlinearity is not covered by our analysis and deserves to be investigated more specifically and carefully.  
\revis{It seems to be very challenging even in the two-dimensional situation.}

The paper is organized as follows. In the following section, we formulate the general assumptions and study the state
system, proving the existence of a unique solution. We also show the uniform separation property for the solution component
$\varphi$ and the local Lipschitz continuity of the control-to-state operator. In Section 3, we then employ the implicit 
function theorem to prove that  the control-to-state operator is twice continuously Fr\'echet differentiable between appropriate 
Banach spaces. Moreover, local Lipschitz properties are shown for the first and second derivatives. In Section 4, the main part
of this paper, we investigate the control problem {\bf (CP)} with sparsity. Besides analyzing the associated adjoint problem, we
derive the first-order necessary optimality conditions. The final section then brings the derivation of the announced second-order 
sufficient optimality conditions for controls that are locally optimal in the sense of $L^2(Q)$.

Prior to this, let us fix some notation.
For any Banach space $X$, we denote by \,$\|\,\cdot\,\|_X$, $X^*$, and $\langle\, \cdot\, , \,\cdot\, \rangle_X$,  
the corresponding norm, its dual space, and  the related duality pairing between $X^*$ and~$X$. 
For two Banach spaces $X$ and $Y$ that are both continuously embedded in some topological vector space~$Z$, we consider the linear space
$X\cap Y$ that becomes a Banach space if equipped with its natural norm $\norma v_{X\cap Y}:=\norma v_X+\norma v_Y\,$  for $v\in X\cap Y$.
The standard Lebesgue and Sobolev spaces defined on $\Omega$ are, for
$1\le p\le\infty$ and $\,m \in \enne \cup \{0\}$, denoted by $L^p(\Omega)$ and $W^{m,p}(\Omega)$, respectively.  
If $p=2$, they become Hilbert spaces, and we use the usual notation $H^m(\Omega):= W^{m,2}(\Omega)$. 
For convenience, we also~set
\begin{align*}
  & H := \Ldue , \quad V := \Huno, \quad W:=\bigl\{v\in \Hdue : \ \dn v =0 \, \hbox{ on }\, \Gamma \bigr\},
  \end{align*}
and we denote by $(\,\cdot\,,\,\cdot\,)_H$ the natural inner product in $H$. 
As usual, $H$ is identified with a subspace of the dual spaces $\VD$ according to the identity
\begin{align*}
	\langle u,v\rangle_V =(u,v)_H
	\quad\mbox{for every $u\in H$ and $v\in V$}.
\end{align*}
We then have the Hilbert triplet $V \subset H \subset \VD$ with dense and compact embeddings.

We close this section by introducing a convention concerning the constants used in estimates within this paper: we denote by $\,C\,$ any 
positive constant that depends only on the given data occurring in the state system and in the cost functional, as well as 
on a constant that bounds the $L^\infty(Q)$--norms of the elements of $\Uad$. The actual value of 
such generic constants $\,C\,$ 	may 
change from formula to formula or even within formulas. Finally, the notation $C_\delta$ indicates a positive constant that
additionally depends on the quantity $\delta$.   


\section{General assumptions and the state system}
\setcounter{equation}{0}

In this section, we formulate the general assumptions for the data of the state system \State, and we collect some known results for
the state system. Throughout this paper, we make the following assumptions:
\begin{description}
\item[(A1)] \,\,$f=f_1+f_2$, where $f_1 :\erre\to [0,+\infty]$ is convex and lower semicontinuous with $f_1(0)=0$, and $f_2:\erre\to\erre$ has 
a Lipschitz continuous first derivative $f_2'$ on $\erre$.
Moreover, we require that $f_1 \in C^5(-1,1)$ and $f_2 \in C^5[-1,1]$, and we assume that
\begin{align}
\label{limf1}
&\lim_{r\searrow -1} \,f_1'(r) =-\infty\,,\quad
\lim_{r\nearrow 1} \,f_1'(r)=+\infty\,.
\end{align}
\item[(A2)] \,\,$\tau>0$,  $\gamma \in L^\infty (\Omega)$, and there exists some $\gamma_0 > 0 $ such that $\gamma\geq \gamma_0\,  $ a.e. in $\Omega$.
Moreover, $w_0\in H$, $\phi_0 \in W $, and it holds that
\begin{align}
\label{ini1}
&-1\,<\, \min_{x\in\overline\Omega}\, \phi_0(x), \quad \,\,\, \max_{x\in\overline\Omega}\,\phi_0(x)<1\,.
\end{align}
\item[(A3)] \,\,$R>0$  is a fixed constant such that 
\Beq
\label{defUR}
\Uad\subset   \CU_R :=\{u\in L^\infty(Q): \,\,\|u\|_{L^\infty(Q)}\,<R\}.
\Eeq  
\end{description} 

\vspace*{1mm}
\Brem
\label{Rem1}
From the condition~{\bf (A1)} (cf.~\eqref{limf1}, in particular) it follows that the derivative $f_1'$ is just defined in $(-1,1)$
 and gives rise to a maximal monotone operator in $\erre \times \erre$.
Note that {\bf (A1)} is fulfilled if $\,f\,$ is given by the logarithmic potential $f_{\rm log}$ in \eqref{flog},
where $f_1(r) = c_1\bigl((1+r)\,\ln(1+r)+(1-r)\,\ln(1-r)\bigr)$ for $r\in (-1,1)$ and $f_2(r) = -c_2r^2$ for $r\in \erre$ in that case. 
The condition $\phi_0 \in W $ implies that $\phi_0$ is uniformly bounded and continuous on $\overline{\Omega}$, so that
\eqref{ini1} yields that $\phi_0 $ is strictly separated from the values $-1$ and $1$ associated with the pure phases.  
Finally, the condition {\bf (A3)} just fixes once and for all a bounded open subset of the control space $\LiQ$ that contains $\Uad$. 
\Erem

A consequence of {\bf (A2)} is that the mean value of $\phiz$,
\Beq
  m_0 := \frac 1 {|\Omega|}\iO \phiz \,,\hbox{ belongs to the interior of the domain $(-1,1) $ of $f_1' $}.
  \label{meanz}
\Eeq
In the following, we use the notation $\overline  v$ to denote the mean value of a generic function $v\in\Luno$.
More generally, we set
\Beq
  \overline v : = \frac 1 {|\Omega|} \, \langle v , 1 \rangle_V
  \quad \hbox{for every $v\in V^*$},
  \label{defmean}
\Eeq
noting that the constant function 1 is an element of $V$. Clearly, $\overline  v$ is the usual mean value of $v$ if $v\in H$.

Next, we specify our notion of solution: for any given $u \in \L2H$, the triplet $(\phi,\mu,w)$ is said to be a solution
to \State\  if 
\begin{align*}
& \phi \in  H^1(0,T;H)\cap C^0([0,T];V) \cap L^2(0,T;W), \\
&-1 < \phi(x,t) < 1\quad \hbox{for a.e. } (x,t)\in Q,\\
& \mu \in \L2W , \quad w \in \H1 H,\\
&(\phi,\mu,w) \ \hbox{ solves \State,}
\end{align*}
so that, in particular, 
\begin{align}
\label{ssvar1}
& \int_\Omega \partial_t \phi (t) v +\int_{\Omega}\nabla \mu (t)\cdot\nabla v =0 
\quad \hbox{for a.e. $t\in (0,T)$ and every $v\in V$,}
\\ 
\label{ssvar2} 
& \tau \int_\Omega \partial_t \phi (t) v +\int_{\Omega}\nabla \phi (t)\cdot\nabla v +\int_{\Omega} f'( \phi (t)) v 
\non \\
&\quad{}
= \int_{\Omega} ( \mu (t)+ w(t) ) v
\quad \hbox{for a.e. $t\in (0,T)$ and every $v\in V$,}
\\
\separa
\label{ssvar3} 
&\int_{\Omega}\gamma \, \dt w (t) z +\int_{\Omega }w(t) z = \int_{\Omega }u(t) z 
\quad \hbox{for a.e. $t\in (0,T)$ and every $z\in H$,}
\end{align} 
as well as 
\Beq
\label{ssvar4} 
\phi(0)=\phi_0   \quad \mbox{in }\,V ,\,\quad w(0)= w_0 \quad \mbox{ in }\,H.
\Eeq
Note that the above identities \eqref{ssvar1}--\eqref{ssvar3} are variational formulations of  \eqref{ss1}--\eqref{ss3}, 
where the first two are obtained with the contribution of the boundary conditions~\eqref{ss4}. Let us emphasize that, 
by this definition, \eqref{ss2} actually holds and, by comparison of terms, it turns out that $f'(\phi)\in L^2(0,T;H)$.
In addition, \eqref{ssvar4} is another way of writing the initial conditions~\eqref{ss5}. 
\revis{By the above notion of solution, it is also clear that $(\phi,\mu,w)$ is a strong solution in the sense of \eqref{ss1}--\eqref{ss5}.}

Let us also remark that, 
thanks to the linear equation~\eqref{ss3} and the second initial condition in \eqref{ssvar3}, $w$ can be explicitly written in terms of $u$ 
by means of the variation of constants formula
\Beq
\label{ssvar5} 
w(x,t)= w_0(x) \exp (-t/\gamma(x)) + \int_0^t \exp(-(t-s)/\gamma(x) ) u(x,s) ds , \quad \hbox{a.e.}\  (x,t) \in Q.
\Eeq

We are going to prove the existence of a (smoother) solution and, in the case when $u\in\CU_R$, the separation property. 
\Bthm
\label{Teo1}
Suppose that the conditions {\bf (A1)}--{\bf (A3)} are fulfilled. Then the state system \State\ has for any 
$u\in \L2H $ a unique solution $(\phi,\mu,w)$ with the regularity
\begin{align}
  & \phi \in \W{1,\infty}H \cap \H1V \cap \L\infty W \,  ,
  \label{regphibis}
  \\
  & \phi \in  C^0(\overline Q) \quad \hbox{and } \,\, -1<\phi <1 \quad \hbox{in }\,  Q \, ,
  \label{regphiter}
  \\
  & \mu \in \L\infty W \cap\L 2{H^3(\Omega)}
  \subset \LQ\infty\, ,
  \label{regmubis}
  \\
  & w\in \H1 H \, . 
  \label{regw}
\end{align}
In addition, there is a constant $K_1>0$, which depends only on $\norma{u}_{\L2H}$ and the data of the state system, such that 
\begin{align}
\label{ssbound1}
&\|\phi\|_{ \W{1,\infty}H \cap \H1V \cap \L\infty W\cap C^0(\overline Q) } \nonumber\\
&+\,\|\mu\|_{ L^\infty(0,T;W)\cap\L 2{H^3(\Omega)}\cap\LQ\infty} \,
+ \, \| w \|_{ \H1 H }       \,\le\,K_1\, , 
\end{align} 
whenever $(\phi,\mu,w)$ is the solution to the state system associated with $u$. 
Moreover, if $w_0\in L^\infty (\Omega)$ and $u\in \CU_R$, then the solution component $w$ satisfies 
\Beq
w \in\W{1,\infty}{L^\infty (\Omega)}  \subset \LQ\infty\,,
  \label{regwbis}
\Eeq
and a uniform strict separation
property is fulfilled: there are constants $r_-,\,  r_+$, which depend only on $\,R\,$ and the data of the state system, such that
\Beq
\label{separation}
-1 < r_- \le \phi(x,t) \le r_+ < 1  \quad\mbox{for every }(x,t)\in \overline Q, 
\Eeq
whenever $\phi$ is first component of the solution $(\phi,\mu,w)$ to the state system related 
to some~$u\in\CU_R$.
\Ethm

\begin{corollary}
\label{Cor1}
Assume {\bf (A1)}--{\bf (A3)}, and let $w_0\in L^\infty (\Omega)$. Then, for all 
$u\in  \CU_R$,  the corresponding  solution $(\phi,\mu,w)$ of the state system \State\ satisfies
\Beq
\label{ssbound2}
\max_{0\le i\le 5}\, \Bigl(\max_{j=1,2}\,\|f_j^{(i)}(\phi)\|_{C^0(\overline Q)} +
\|f^{(i)}(\phi)\|_{C^0(\overline Q)}  \Bigr)\,\le\,K_2\,,
\Eeq
for some constant $K_2$ depending only on $r_-$, $r_+$, $ f_1$, $f_2$, where $f^{(i)}= f_1^{(i)}+ f_2^{(i)}$ for~$i=0,1,\ldots, 5$. 
\end{corollary}

Notice that the regularity $\varphi\in C^0(\overline Q)$ follows from \eqref{regphibis} and \cite[Sect.~8,~Cor.~4]{Simon},
since the continuous embedding $W\subset C^0(\overline\Omega)$ is compact. The estimate \eqref{ssbound2} is then an 
immediate consequence of \eqref{separation} and assumption {\bf (A1)}. The proof of the above 
theorem, however, is rather long and involved.

\begin{proofteo1}
To begin with, we consider for every $\varepsilon\in (0,1)$ the Moreau--Yosida regularization \revis{$\Betaeps$ of $f_1$ (see, e.g., \cite{Barbu, Brezis}),  which is defined by
\begin{align*}
	& \Betaeps(r)
	:=\inf_{s \in \mathbb{R}}\left\{ \frac{1}{2\varepsilon } |r-s|^2
	+f_1(s) \right\} 
	= 
	\frac{1}{2\varepsilon } 
	\bigl| r-R_\varepsilon  (r) \bigr|^2+f_1 (R_\varepsilon (r) )
	= \int_{0}^{r} \betaeps (s)ds 
\end{align*}
for all $r\in\erre$. We point out that  $\betaeps :\mathbb{R} \to \mathbb{R}$ \juerg{and} the associated resolvent operator
$R_\varepsilon :\mathbb{R} \to \mathbb{R}$ are represented by 
\begin{align*}
	& \betaeps (r)
	:= \frac{1}{\varepsilon } ( r-R_\varepsilon (r) ), 
	\quad 
	R_\varepsilon (r) 
	:=(I+\varepsilon  f_1' )^{-1} (r)\juerg{, \quad\mbox{for all \,$r\in\erre$,}}
\end{align*}
\juerg{with} $I$ denoting here the identity operator. Note that the derivative $\betaeps$, defined in $\erre$,
turns out to be the Yosida approximation of the maximal monotone graph induced by 
$f_1' $, that has effective domain $(-1,1)$. Now, in view of {\bf (A1)} (cf. also Remark~\ref{Rem1}), 
$\betaeps$ and $\Betaeps$ satisfy  the following properties} (see, e.g., \cite[pp.~28 and~39]{Brezis}): 
\begin{align}
  & \hbox{$\betaeps :\erre \to \erre $ is monotone and \Lip\ continuous} 
	\nonumber\\
	&\quad \hbox{with Lipschitz constant } 1/\eps, \mbox{ and it holds} \,\betaeps(0)=0,
  \label{monbetaeps}
  \\
  & |\betaeps(r)| \leq |f_1^\prime (r)|
  \quad \hbox{for every $r\in (-1,1)$},
  \label{disugbetaeps}
  \\
  & \revis{{}0 \leq \Betaeps(r) \leq f_1(r){}}
  \quad \hbox{for every $r\in\erre$}.
  \label{disugBetaeps}
\end{align}
\Accorpa\Propbetaeps monbetaeps disugBetaeps
Then, consider the problem of finding $(\phieps ,\mueps , w) $ satisfying \eqref{ssvar5} as well~as
\begin{align}
  & \iO \dt\phieps(t) \, v 
  + \iO \nabla\mueps(t) \cdot \nabla v
  = 0 
  \quad \hbox{for a.e. $t\in (0,T)$ and every $v\in V$,}
  \qquad
  \label{primaeps}
  \\
  \separa
  & \tau \iO \dt\phieps(t) \, v
  + \iO \nabla\phieps(t) \cdot \nabla v
  + \iO \bigl( \betaeps(\phieps(t)) + f'_2 (\phieps(t)) \bigr) v
  \non
  \\
  & \quad {}
  = \iO \bigl( \mueps(t) + w(t) \bigr) v
  \quad \hbox{for a.e. $t\in (0,T)$ and every $v\in V$,}
  \label{secondaeps}
  \\
  \separa
  & \phieps(0) = \phi_0 \quad\hbox{a.e. in $\Omega$.}
  \label{cauchyeps}
\end{align}
\Accorpa\Pbleps primaeps cauchyeps
Note that \eqref{primaeps}--\eqref{cauchyeps} is a well-known viscous Cahn--Hilliard system that has received a lot of attention in the recent literature. In addition, here 
the nonlinearies acting in \eqref{secondaeps} are even Lipschitz continuous. For the existence and uniqueness of  
a solution $(\phieps ,\mueps) $ to \eqref{primaeps}--\eqref{cauchyeps},
we may refer to \cite[Thm.~4.1]{CGM}, where a Faedo--Galerkin scheme has been employed for the proof. Observe that the weak solution
offered by~\cite[Thm.~4.1]{CGM} is a variational solution with the regularity
$\phieps \in H^1(0,T;V^*)\cap C^0([0,T];H) \cap L^2(0,T;V)$ and $ \mueps \in \L2V$; however, since 
$$ - \betaeps(\phieps) - f'_2 (\phieps)
 + \mueps + w \in \L2V + \H1H \quad\hbox{and} \quad  \phi_0 \in W , 
$$ 
it follows from classical parabolic regularity theory (see, e.g., \cite{Lions}) that 
\begin{equation}
\phieps \in  \H1V \cap \L\infty W \quad\hbox{and} \quad  \mueps \in \L2W, 
\label{p-1}
\end{equation}
at least. Moreover, the pointwise equations and conditions
\begin{align}
\label{sseps1}
&\dt \phieps-\Delta \mueps = 0 &&\mbox{a.e. in }\,Q,
\\
\label{sseps2}
&\tau \dt \phieps -\Delta \phieps + \betaeps(\phieps) + f'_2 (\phieps)   = \mueps + w &&\mbox{a.e. in }\,Q,
\\
\label{sseps3}
&\dn \mueps = \dn\phieps = 0 &&\mbox{a.e. on }\,\Sigma,
\end{align}
are valid. However, in the sequel the reader can realize how to directly 
obtain the regularities in \eqref{p-1} and even more. Indeed, we
are now going to recover a number of a priori estimates, where $C>0$ denotes constants
 that are independent of $\eps\in (0,1)$. 

\begin{step}
First estimate

Take $v=1/|\Omega| $ in \eqref{primaeps} and integrate with respect to time using \eqref{cauchyeps}. Recalling \eqref{meanz}
and \eqref{defmean}, we obtain the mean value conservation property
\Beq 
\juerg{\overline{\phieps(t)}} = \overline{\phi_0} = m_0 \quad \hbox{for all } t\in [0,T].
\label{pier1}
\Eeq
Now, let us make a preliminary remark for $w$. As $0\leq \exp (-t/\gamma(x))\leq 1$ for all $t\geq0$,  
 it is not difficult to deduce from \eqref{ssvar5} and \eqref{ss3} that 
\begin{align}
&\label{pier2}
\norma{w}_{\L\infty H} \leq \norma{w_0}_{H} + \sqrt{T} \, \norma{u}_{\L2H},
\\
&\label{pier3}
\norma{\dt w}_{\L2H} \leq \frac1{\gamma_0} \Bigl(\norma{u}_{\L2H} + \sqrt{T}\, \norma{w_0}_{H} + T\norma{u}_{\L2H}\Bigr).
\end{align} 
Next, we choose $v= \mueps (t) $ in \eqref{primaeps} and $v=\dt \phieps (t) $ in \eqref{secondaeps}, add the resulting equalities 
and integrate with respect to t. Noting that a cancellation occurs, we obtain 
\begin{align}
&\iQt |\nabla\mueps|^2 +\tau \iQt |\dt \phieps|^2 + \frac12 \iO |\nabla \phieps(t)|^2 + \iO \Betaeps (\phieps (t) ) \non
\\
&=  \frac12 \iO |\nabla \phiz|^2 + \iO \betaeps (\phiz ) + \iQt (w -  f'_2 (\phieps)) \dt \phieps. 
\label{pier4}
\end{align}
The first two terms on the \rhs\ are under control due to {\bf (A1)}, {\bf (A2)} and \revis{\eqref{disugbetaeps}}. For the third term,
 we have that
 $$  \iQt (w -  f'_2 (\phieps)) \dt \phieps \leq \frac \tau 2 \iQt |\dt \phieps|^2  + C  \iQt \bigl( 1+ |\phieps|^2 + |w|^2 \bigr), $$
thanks to Young's inequality and the Lipschitz continuity of $f'_2$. 

Now observe that, thanks to \eqref{pier1} and the Poincar\'e--Wirtinger inequality, there is some constant $c_1>0$ depending only on $\Omega$ such 
that 
$$
\frac12 \iO |\nabla \phieps(t)|^2\,\ge\,c_1 \norma{\phieps(t)- m_0}^2_V.
$$  
Therefore, by combining the above estimate with \eqref{pier4}, we can infer that
\begin{align}
&\iQt |\nabla\mueps|^2 +\frac \tau 2 \iQt |\dt \phieps|^2 + c_1 \norma{\phieps(t)- m_0}^2_V + \iO \Betaeps (\phieps (t) ) \non
\\
&\leq \frac12 \iO |\nabla \phiz|^2 + \iO \betaeps (\phiz ) + C\norma{w_0}^2_{H} 
+ C  \iQt\bigl( 1+ |\phieps(s)-m_0|^2 \bigr) + C \norma{u}^2_{\L2H}\,,
\non
\end{align}
whence, in view of \eqref{pier2} and Gronwall's lemma, it is straightforward to arrive at the estimate 
\Beq
\label{pier5}
\norma{\nabla\mueps}_{\revis{\L2 {\revis{H}}^3}} + \norma{\phieps}_{\H1H\cap \L\infty V} +\norma{ \Betaeps (\phieps ) }_{\L\infty{L^1(\Omega)}} \leq C.
\Eeq
\end{step}

\begin{step}
Second estimate

We take $v=\phieps(t) - m_0 $ in \eqref{secondaeps}, in order to exploit
the following argument, which owes to \cite[Appendix, Prop.~A.1]{MiZe} (see also \cite[p.~908]{GiMiSchi} for a detailed proof) 
and was used in several papers: by virtue of~\eqref{meanz}, we have, with some $\delta_0>0$ depending only on $f_1'$ and~$m_0$, that
\Beq
  \betaeps({r}) (r-m_0)
  \geq \delta_0 |\betaeps(r)| - \delta_0^{-1}
  \quad \hbox{for every $r\in\erre$ and every $\eps\in(0,1)$}.
  \label{trickMZ}
\Eeq
We argue for fixed $t$ and avoid time integration, where, for simplicity, we do not write the time $t$ for a while.
We have, almost everywhere in $(0,T)$,
\begin{align}
  & {\delta_0}\iO |\betaeps(\phieps)| - \delta_0^{-1} \, |\Omega|
  \leq \iO \nabla\phieps \cdot \nabla(\phieps-m_0)
  + \iO \betaeps(\phieps) (\phieps-m_0)
  \non
  \\
  & = \iO \mueps (\phieps-m_0)
  + \iO \bigl(
    w
    - \tau \dt\phieps
    - f_2'(\phieps)
  \bigr) (\phieps-m_0) .
  \label{perpt}
\end{align}
We recall that $\overline{\phieps-m_0}=0$ a.e. in $(0,T)$, thus we can take advantage 
of that and apply the Poincar\'e--Wirtinger inequality {to 
$\mueps-\overline\mueps$. In fact, using \eqref{pier5} as well,} we have
\begin{align}
  & \iO \mueps (\phieps-m_0)
  = \iO (\mueps-\overline\mueps) (\phieps-m_0)
  \leq C \, \norma{\nabla\mueps}_{H^3} \, \norma{\phieps-m_0}_H 
  \leq C \, \norma{\nabla\mueps}_{H^3}  \,.
  \non
\end{align}
For the remaining terms on the \rhs\ of \eqref{perpt},
we use the Schwarz inequality, the Lipschitz continuity of $f'_2$, and the bounds available from~\eqref{pier2} and \eqref{pier5}, to obtain that
\begin{align}
  & \norma{\betaeps(\phieps)}_1
  \leq C\, \bigl(
    \normaH{\nabla\mueps} 
    + \normaH{\dt\phieps}
    + 1
  \bigr) \,\mbox{ a.e. in }\,(0,T).
  \non
\end{align}
At this point, by taking $v=1/|\Omega|$ in \eqref{secondaeps}, using the inequality just obtained,
and estimating the other $L^1$-norms by the corresponding $H$-norms, we deduce that
\begin{align}
  |\overline\mueps|
  \leq C \, \bigl(
    \normaH{\nabla\mueps} 
    + \normaH{\dt\phieps}
    + 1
  \bigr) \quad \hbox{a.e. in $(0,T)$.}
  \label{point-wise}
\end{align}
Then, by \eqref{point-wise} and \eqref{pier5} we find that 
$
\norma{\mueps}_{\L2 {V}}  \leq C.
$
Moreover, using again the boundedness of $\dt \phieps$ in $\L2 H$ along with elliptic regularity theory, 
we additionally recover from \eqref{sseps1} and \eqref{sseps3}  the estimate 
\Beq
\label{pier6}
\norma{\mueps}_{\L2 {W}}  \leq C.
\Eeq

\end{step}

\begin{step}
Third estimate

We take $v= \betaeps (\phieps (t)) $ as test function in \eqref{secondaeps} and do not integrate 
with respect to time, obtaining for a.e. $t\in (0,T)$ that
\begin{align}
  & 
  \iO \nabla\phieps(t) \cdot \nabla \betaeps (\phieps (t))
  + \iO |\betaeps (\phieps (t))|^2 
  \non
  \\
  & = \iO \bigl( \mueps 
    + w
    - \tau \dt\phieps
    - f_2'(\phieps)
  \bigr)(t)  \, \betaeps (\phieps (t)).
  \label{pier7}
\end{align}
Note that the first term on the \lhs\ is nonnegative,
 while the \rhs\ can be easily treated using Young's inequality and taking advantage of \eqref{pier2} and \eqref{pier5}, to deduce that
\begin{align}
  & \frac 12 \norma{\betaeps(\phieps(t))}^2_H
  \leq C_R\, \Bigl(
    \norma{\mueps(t) }^2_H 
    + \normaH{\dt\phieps (t) }^2
    + 1
  \Bigr) \quad \hbox{for a.e. $t\in (0,T)$}.
  \label{pier8}
\end{align}
Moreover, since it follows from \eqref{sseps2}  that 
$$ 
-\Delta \phieps(t) = -  \betaeps(\phieps(t) ) + \bigl(\mueps + w - \tau \dt \phieps + f'_2 (\phieps)\bigr)(t) \quad \hbox{a.e. in $\Omega$}, 
$$ 
we can infer from \eqref{pier8} that 
\begin{align}
  & 
   \norma{\Delta \phieps (t))}^2_H
  \leq C\, \Bigl(
    \norma{\mueps(t) }^2_H
    + \normaH{\dt\phieps (t) }^2
    + 1
  \Bigr) \quad \hbox{for a.e. $t\in (0,T)$}.
  \label{pier9}
\end{align}
Hence, owing to \eqref{sseps3}, \eqref{pier5}, \eqref{pier6} and elliptic regularity theory, we conclude that 
\Beq
\label{pier10}
\norma{\betaeps(\phieps)}_{\L2 {H}} +  \norma{\phieps}_{\L2 {W}}  \leq C.
\Eeq
\end{step}

\begin{step}
Fourth estimate

The subsequent estimate should be rigorously reproduced on some regularized version of 
\eqref{sseps1}--\eqref{sseps3} with \eqref{cauchyeps}; for instance, one can use a time discetization procedure. 
However, for the sake of brevity, let us argue directly on \eqref{sseps1}--\eqref{sseps3}. 
First, we write \eqref{sseps1}--\eqref{sseps2} at the initial time $t=0$ and deduce from 
\eqref{sseps2}, \eqref{cauchyeps} and \eqref{ssvar4} that
\Beq 
\tau \dt \phieps (0) =  \Delta \phiz - \betaeps(\phiz) - f'_2 (\phiz) + \mueps (0) + w_0 , \label{pier11}
\Eeq
whence, replacing $\dt \phieps (0)$ in \eqref{sseps1}, we obain the elliptic equation
\Beq 
\mueps (0)- \tau \Delta \mueps (0) =  - \Delta \phiz + \betaeps(\phiz) + f'_2 (\phiz) - w_0,
\label{pier11bis}
\Eeq
where the \rhs\  is bounded in $H$ due to {\bf (A2)}, \eqref{limf1}, and  \eqref{disugbetaeps}.
Then, from the homogeneous boundary condition $\dn \mueps (0)=0$ on $\Gamma$ (see \eqref{sseps3}) and the elliptic
well-posedness and regularity theory,
it turns out that there exists a unique solution $\mueps (0)$ to \eqref{pier11bis} satisfying
\Beq
\label{pier12}
\norma{\mueps (0)}_{W}  \leq C.
\Eeq
Moreover, coming back to \eqref{pier11}, we also recover that 
\Beq
\label{pier13}
\norma{\dt \phieps (0)}_{H}  \leq C.
\Eeq
The next (formal) computation is performed directly on the variational formulation \linebreak 
\eqref{primaeps}--\eqref{secondaeps} of \eqref{sseps1}--\eqref{sseps3}. It
consists in differentiating \eqref{secondaeps} with respect to $t$ and then 
taking $v= \dt \phieps$. On the other hand, we choose   $v= \dt \mueps $ in \eqref{primaeps}
and add the result to the previous equality. Note that a cancellation occurs. Then, we integrate over $(0,t)$
and obtain  
\begin{align}
&\frac12 \iO |\nabla\mueps (t) |^2 + \frac\tau 2 \iO |\dt \phieps(t) |^2 + \iQt |\nabla \dt \phieps|^2 + \iQt \dt 
(f_{1,\varepsilon}' (\phieps )) \dt \phieps \non
\\
&\leq  \frac12 \iO |\nabla\mueps (0) |^2 + \frac\tau 2 \iO |\dt \phieps(0) |^2
+ \iQt \dt (w -  f'_2 (\phieps)) \dt \phieps. 
\non
\end{align}
In view of \eqref{monbetaeps}, the fourth term on the \lhs\ is nonnegative. 
On the \rhs , we invoke \eqref{pier12} and \eqref{pier13} for the first two terms and easily control the third one by 
using \eqref{pier3}, the Lipschitz continuity of $f'_2$, and the boundedness of 
$\norma{\dt \phieps}_{\L2H}^2 $ established in \eqref{pier5}. Hence, we easily conclude that 
\Beq
\label{pier14}
\norma{\nabla\mueps}_{\revis{\L\infty {H}^3}} + \norma{\phieps}_{\W{1,\infty}H\cap \H1 V} \leq C.
\Eeq
\end{step} 

\begin{step}
Conclusion of the existence proof

Next, we return to the bound \eqref{point-wise} for the mean value of $\mueps$, 
observing that by virtue of \eqref{pier14} the \rhs\ of \eqref{point-wise} is now bounded in $L^\infty (0,T)$. 
Hence, we have that 
$
\norma{\mueps}_{\L\infty {V}}  \leq C,
$
and, in view of the boundedness of $\dt \phieps$ in $\L\infty H \cap \L 2 V$ from \eqref{pier14}, 
we can exploit the equation \eqref{sseps1}, the boundary condition in \eqref{sseps3}, and the elliptic regularity theory, to arrive at
\Beq
\label{pier15}
\norma{\mueps}_{\L\infty {W}\cap \L2 {H^3(\Omega)}}  \leq C.
\Eeq
Similarly, we can recall \eqref{pier8} and \eqref{pier9} and, with the help of \eqref{pier14} and \eqref{pier15}, improve the estimate \eqref{pier10}, deducing that 
\Beq
\label{pier16}
\norma{\betaeps(\phieps)}_{\L\infty {H}} +  \norma{\phieps}_{\L\infty {W}}  \leq C.
\Eeq

At this point, we can perform the passage to the limit as $\eps\searrow0$. 
In view of the estimates \eqref{pier14}--\eqref{pier16}, {which are independent of $\eps$,} by weak and weak-star compactness it turns out that there are $\mu,\phi$ and $\zeta$ such that 
\begin{align}
&  \phieps \to \phi \quad \hbox{weakly star in } \   \W{1,\infty}H\cap \H1 V \cap \L\infty W , \label{pier16-1}
\\
& \mueps \to \mu \quad \hbox{weakly star in } \  \L\infty {W}\cap \L2 {H^3(\Omega)},  \label{pier16-2}
\\
& \betaeps(\phieps) \to \zeta \quad \hbox{weakly star in } \   \L\infty {H}, \label{pier16-3}
\end{align} 
as $\eps \searrow 0$, possibly along a subsequence. By virtue of \eqref{pier16-1}--\eqref{pier16-3} and the Aubin--Lions--Simon
 lemma (see, e.g., \cite[Sect.~8, Cor.~4]{Simon}, as $W \subset C^0 (\overline\Omega)$ with compact embedding), 
we deduce in particular that $\phieps\to \phi$ strongly in $C^0 (\overline Q) $. The same strong convergence holds 
for $f'_2 (\phieps) \to f'_2(\phi) $, while the identification $\zeta = f'_1 (\phi)$ results 
(first as an inclusion) as a consequence of \eqref{pier16-3} and the maximal monotonicity of $f'_1$ 
(see {\bf (A1)} and Remark~\ref{Rem1}), since we can apply, e.g., \cite[Lemma~2.3, p.~38]{Barbu}.
Then, we can pass to the limit in the variational equalities  \eqref{primaeps}, \eqref{secondaeps} and obtain \eqref{ssvar1}, \eqref{ssvar2}. Also, the initial condition 
\eqref{cauchyeps} extends to the limit $\phi$. Finally, we have found a complete solution $(\phi,\mu, w)$ to \State, solving then 
\eqref{ssvar1}--\eqref{ssvar4}, possessing the full regularity expressed in \eqref{regphibis}--\eqref{regw} and, due to 
\eqref{pier16-1}--\eqref{pier16-3} and the weak star lower semicontinuity of norms, satisfying the bound \eqref{ssbound1}. 
The uniqueness of this solution $(\phi,\mu, w)$ will follow as a consequence of the subsequent continuous dependence result. 

\end{step}

\begin{step}
Separation property

Now assume, in addition to {\bf (A1)--(A3)}, that $w_0\in L^\infty (\Omega)$ and $u\in \CU_R$. We aim at verifying 
the separation property~\eqref{separation}. First, note that \eqref{regwbis} is a direct consequence of \eqref{ssvar5} and {\bf (A3)}. 
Moreover, the equation \eqref{ss2} holds for the limit 
functions with the datum $f'= f_1' + f_2'$ as in {\bf (A1)}. Therefore, we rewrite \eqref{ss2} in the form
\Beq
\label{pier17} 
\tau \dt \phi -\Delta \phi + f'_1 (\phi) = \mu + w - f'_2(\phi)  \quad\mbox{a.e. in }\,Q,
\Eeq
noting that the \rhs\ is bounded in  $L^\infty(Q)$ (cf. \eqref{regphiter} and \eqref{regmubis}). In fact,
there exists a positive constant $c_*$, independent of the choice of $u\in {\cal U}_R$, such that 
\begin{equation} 
\label{pier18}
\norma{\mu + w - f'_2(\phi) }_{L^\infty(Q)}\leq c_* .
\end{equation}
Moreover, the condition \eqref{ini1} for the initial value $\phi_0$ and the assumption~\eqref{limf1} entail the existence of 
some constants $r_-$ and $r_+$ such that $-1 < r_-\leq r_+ < 1 $
and
\begin{gather}
	\label{separation_first}
	 r_- {{}\leq{}} \min_{x\in \overline\Omega} \phi_0(x), \quad r_+ {{}\geq{}} \max_{x\in \overline\Omega}\, \phi_0(x),
	\\
	f'_1 (r)
	 + c_* \leq 0 
	\quad \forall r \in (-1,r_-), \quad 
	f'_1(r)
	 - c_* \geq 0 
	\quad \forall r \in (r_+,1).
	\label{separation_second}
\end{gather}
Now, let us test \eqref{pier17}, witten at time $t\in(0,T)$, by $v (t) = (\phi(t) - r_+)^+ - (\phi(t) - r_-)^-$, 
where  $(\,\cdot\,)^+$ and $(\,\cdot\,)^- $ denote the
standard positive and negative parts, respectively.
Then, {we} integrate with respect to $t$. Observe that  $v(0)=0$ in view of  \eqref{separation_first}. 
Using integration by parts and \eqref{pier18}, it is straightforward to deduce that
\begin{align*}
	&\frac\tau 2 \norma{v(t)}^2 +  \iQt |\nabla v|^2 \\
	&{}=  \iint_{Q_t \cap \{\phi > r_+\}} \bigl((\mu + w - f'_2(\phi))- f'_1 (\phi)\bigr) (\phi - r_+)\\
	&\quad{}
	+ \iint_{Q_t \cap \{\phi < r_-\}} \bigl(f'_1 (\phi) - (\mu + w - f'_2(\phi)\bigr) (r_- - \phi ) 
 \\
	& \leq{}\iint_{Q_t \cap \{\phi > r_+\}} (c_* - f'_1(\phi)) (\phi - r_+) + \iint_{Q_t \cap \{\phi < r_-\}}
	(f'_1(\phi) + c_*) (r_-- \phi ) .
\end{align*}
Note that the quantity on the last line above is nonpositive due to 
\eqref{separation_second}, so that $v =0$ almost everywhere,
which in turn implies~that 
\begin{align}
	\non
	r_- \leq {\phi} \leq  r_+ \quad \hbox{a.e. in } Q.
\end{align}
Thus, since $\phi\in C^0(\overline Q)$, the separation property \eqref{separation} holds true,
which completes the proof of the assertion.
\end{step}
\end{proofteo1}

Next, we state a continuous dependence result that, in particular, guarantees the uniqueness of the solution provided by Theorem~\ref{Teo1}. 

\Bthm
\label{Teo2}
Suppose that the conditions {\bf (A1)}--{\bf (A3)} are fulfilled. 
If $u_i\in\L2H $, $i=1,2$, are given  
and $(\phi_i,\mu_i,w_i)$, $i=1,2$, are the corresponding solutions to \State, then
\begin{align}
  &\norma{\phi_1-\phi_2}_{\C0H\cap\L2V}
   + \norma{w_1-w_2}_{\H1H}
  \leq K_3 \, \norma{u_1-u_2}_{\L2H},
  \label{contdep1}
\end{align}
for some constant $K_3 $ depending only on $\tau$, $\gamma_0$, $T$ and the Lipschitz constant of $f'_2.$
If, in addition, $w_0\in L^\infty (\Omega)$ and $u_i\in \CU_R$, $i=1,2$, then we have
the further estimate
\begin{align}
  &\norma{\phi_1-\phi_2}_{\H1H\cap\C0V\cap\L2W }
    + \norma{\mu_1-\mu_2}_{\L2W} 
		\non\\
	&\quad	+ \norma{w_1-w_2}_{\H1H} 
  \leq K_4 \, \norma{u_1-u_2}_{\L2H},
  \label{contdep2}
\end{align}
with a constant $K_4$ that depends only on $K_2$, $\tau$, $\gamma_0$,  
$\Omega$, and $T$.
\Ethm

\Bdim
Let us set, for convenience,
\Bsist
  && u := u_1 - u_2 \,, \quad
 \phi := \phi_1 - \phi_2 \,, \quad
  \mu := \mu_1 - \mu_2 \,, \quad
  w := w_1 - w_2\,.
  \non
\Esist
Then $\phi (0) = w(0) = 0 $ in $H$ by \eqref{ssvar4}. 
Next, we can write \eqref{ssvar1}--\eqref{ssvar3} for $(\phi_i,\mu_i,w_i)$, $i=1,2$, and take the differences, obtaining
\begin{align}
\label{diff1}
& \int_\Omega \partial_t \phi\, v +\int_{\Omega}\nabla \mu \cdot\nabla v =0 
\quad \hbox{for every $v\in V$, a.e. in $(0,T)$,}
\\ 
\label{diff2} 
& \tau \int_\Omega \partial_t \phi\,  v +\int_{\Omega}\nabla \phi\cdot\nabla v 
+\int_{\Omega} (f'_1 ( \phi_1) -  f'_1 ( \phi_2) ) v 
- \int_{\Omega} \mu v
\non \\
&\quad{}
= \int_{\Omega} ( w - f'_2 ( \phi_1) +  f'_2 ( \phi_2)) v
\quad \hbox{for every $v\in V$, a.e. in $(0,T)$,}
\\
\separa
\label{diff3} 
&\int_{\Omega}\gamma \, \dt w \, z +\int_{\Omega }w z = \int_{\Omega }u z 
\quad \hbox{for every $z\in H$, a.e. in $(0,T)$.}
\end{align} 
Now, we integrate \eqref{diff1} over $(0,t)$, for $t\in (0,T]$, and take $v= \mu(t)$. At the same time, we insert 
$v = \phi(t)$ in \eqref{diff2} and $z= \dt w (t)$  in \eqref{diff3}. Then we add the three resultants, 
noting that a cancellation occurs, and integrate with respect to time. All this leads to 
\begin{align}
&\frac 12 \iO \Bigl|\nabla\int_0^t \mu(s) ds \Bigr|^2 + \frac \tau 2  \iO |\phi (t) |^2 
+ \iQt |\nabla \phi|^2 
\non \\
&\quad {}+ \iQt \bigl(f'_1 ( \phi_1) -  f'_1 ( \phi_2) \bigr)(\phi_1 - \phi_2)
+ \iQt \gamma |\dt w|^2 +\frac12 \iO |w(t)|^2
\non
\\
&{}=  \iQt \bigl( w \phi  - \bigl(f'_2 ( \phi_1) -  f'_2 ( \phi_2)\bigr) (\phi_1 - \phi_2) + u \, \dt w \bigr) \,,
\label{pier19}
\end{align}
for every $t\in [0,T]$. We point out that the fourth term on the \lhs\ is nonnegative, due to the monotonicity of $f'_1$.
For the \rhs, we infer from the Lipschitz continuity of $f'_2$ and Young's inequality that 
\begin{align}
&  \iQt \bigl( w \phi  - \bigl(f'_2 ( \phi_1) -  f'_2 ( \phi_2)\bigr) (\phi_1 - \phi_2) + u \, \dt w\bigr)
\non\\
&
\leq \frac{\gamma_0}2 \iQt |\dt w|^2 + C \int_0^t \bigl(\normaH{w(s)}^2 + \normaH{\phi(s)}^2
+ \normaH{u(s)}^2\bigr){ds}\,. 
\label{pier20}
\end{align}
Then, as $\gamma \geq \gamma_0$ a.e.~in $\,\Omega\,$ by {\bf (A2)}, we can combine \eqref{pier20} 
with \eqref{pier19} and arrive at \eqref{contdep1} via an application of Gronwall's lemma. 

Now let $w_0\in L^\infty (\Omega)$ and $u_i\in \CU_R$, $i=1,2$. We can then exploit the separation property 
\eqref{separation} and the
global bound \eqref{ssbound2} from Corollary~\ref{Cor1}, which ensures the Lipschitz continuity of $f'=f_1'+f_2'$. 
Now, we test \eqref{diff1} by $\mu$ and 
\eqref{diff2} by $\dt \phi$, then we add and integrate over $(0,t)$, for $t\in (0,T]$. Then, we easily obtain that
\begin{align}
&\iQt |\nabla\mu|^2 +\tau \iQt |\dt \phi|^2 + \frac12 \iO |\nabla \phi(t)|^2  
\leq   \iQt  \bigl( K_2 |\phi| + |w| \bigr) \, |\dt \phi| \non \\
&\quad{}\leq  \frac\tau 2 \iQt |\dt \phi |^2   + C\bigl( \norma{\phi}^2_{L^2(0,t;H)} + \norma{w}^2_{L^2(0,t;H)} \bigr).
\label{pier21}
\end{align}
Taking advantage of \eqref{contdep1}, we then deduce that
\begin{align}
  &\norma{\nabla (\mu_1-\mu_2)}_{\revis{\L2{H}^3}} + \norma{\phi_1-\phi_2}_{\H1H\cap\C0V}
  \leq C \, \norma{u_1-u_2}_{\L2H}. 
  \label{contdep3}
\end{align}

At this point, we can take $v=1/|\Omega|$ in \eqref{diff2} to produce an estimate for the mean value of
$\mu$, since all the other terms are under control. Indeed, it is clear that 
$$
|\overline{\mu (t)}| \leq C \Bigl( \normaH{ \dt \phi(t) } + \normaH{ f' ( \phi_1(t)) -  f' ( \phi_2(t) )} + \normaH{ w (t) } \Bigr) 
$$
for a.e. $t\in (0,T)$; hence it follows from the Lipschitz continuity of $f'$ and from \eqref{contdep1} and \eqref{contdep3} 
that $\norma{\overline{\mu}}_{L^2(0,T)} \leq C\norma{u}_{\L2H}$. Moreover,
by a comparison of terms in \eqref{diff1}, it turns out that $\norma{\Delta \mu }_{\L2H} \leq C\norma{\dt \phi}_{\L2H}$. Therefore, 
first using the Poincar\'e--Wirtinger inequality and then elliptic regularity, we can conclude that 
\begin{align}
  &\norma{\mu_1-\mu_2}_{\L2W} \leq C \, \norma{u_1-u_2}_{\L2H}. 
  \label{contdep4}
\end{align}
Finally, recall that \eqref{diff2} yields 
\Beq
\label{pier22} 
\tau \dt \phi -\Delta \phi = -f'(\phi_1) + f'(\phi_2) + \mu + w  \quad\mbox{a.e. in }\,Q,
\Eeq
as well as $\dn \phi=0$ a.e. on $\Sigma$. Therefore, a comparison of terms in \eqref{pier22} 
leads to the estimate 
$$ \norma{\Delta \phi }_{\L2H} \leq C \Bigl(\norma{\phi}_{\H1H} + \norma{\mu}_{\L2H} + \norma{w}_{\L2H} 
\Bigr), $$
whence,  by virtue of the previous estimates \eqref{contdep1}, \eqref{contdep3}, \eqref{contdep4}, and 
using elliptic regularity theory, we obtain that 
\begin{align}
  &\norma{\phi_1-\phi_2}_{\L2W} \leq C \, \norma{u_1-u_2}_{\L2H}. 
  \label{contdep5}
\end{align}
With	 this, \eqref{contdep2} is completely proved. 
\Edim

\Brem
1.~Note that, by virtue of Theorems~\ref{Teo1} and~\ref{Teo2}, the control-to-state operator $$\CS \revis{{}= (\CS_1, \CS_2 , \CS_3) {}} :u\mapsto \CS(u) \revis{{}= (\CS_1(u), \CS_2(u) , \CS_3(u) ) {}}:= (\phi,\mu, w)$$ is Lipschitz continuous on the set $\CU_R$  as a mapping
between $\L2H$ and the Banach space 
$$\left(\H1H\cap\C0V\cap\L2W \right)\times {\L2W} \times {\H1H}.$$
2.~In view of Theorem~\ref{Teo2}, we emphasize that the continuous dependence estimate 
\eqref{contdep1}, which does not use the bounds \eqref{ssbound1} and \eqref{ssbound2}, is already enough to ensure the uniqueness of the solution $(\phi,\mu, w)$ to \State . Indeed, the uniqueness 
of $\phi$ and $w$ comes immediately from \eqref{contdep1}, while the uniqueness 
of $\mu$ follows from a comparison in \eqref{ss2}, since $f'_1$ is single-valued in its domain $(-1,1)$.
\Erem


\section{Differentiability of the control-to-state operator}
\setcounter{equation}{0}

In this section, we study the differentiability properties of the control-to-state operator $\CS$. In addition to the control
space $\,{\cal U} = \LiQ \,$ defined in \eqref{defU}, we introduce the
Banach spaces
\begin{align}
\label{defX}
&\CX:=\left(H^1(0,T;H)\cap C^0([0,T];V)\cap L^2(0,T;W)\right)\times L^2(0,T;W) \times H^1(0,T;H),\\
\label{defY}
&\CY:=\left(W^{1,\infty}(0,T;H)\cap H^1(0,T;V)\cap L^\infty(0,T;W)\right)\nonumber\\
&\hspace*{11mm} \times \left(L^\infty(0,T;W)\cap L^2(0,T;H^3(\Omega))\right)\times H^1(0,T;H),\\
\label{defZ}
&\CZ:=\left\{(\phi,\mu,w)\in \CY\cap \CU^3: 
\dt\phi-\Delta\mu\in\CU, \right.\nonumber\\
&\hspace*{3cm} \left.\tau\dt\phi-\Delta\phi-\mu-w\in\CU , \,\,\,\gamma\dt w+w\in\CU\right\}, 
\end{align}
where $\CX$ and $\CY$ are endowed with their standard norms and the norm in $\CZ$ is given by
\begin{align}
&\|(\phi,\mu,w)\|_{\CZ}\,=\,\|(\phi,\mu,w)\|_{\CY} + \|(\phi,\mu,w)\|_{\CU^3} + \|\dt\phi-\Delta\mu\|_{\CU} \nonumber\\
&\hspace*{29mm} +\,\|\tau\dt\phi-\Delta\phi-\mu-w\|_{\CU} + \|\gamma \dt w+w\|_{\CU}. 
\end{align}  
We want to show that under the assumptions {\bf (A1)}--{\bf (A3)} and $w_0 \in L^\infty (\Omega)$ 
the operator $\CS$ is  twice continuously Fr\'echet differentiable on $\CU$ 
as a mapping from $\CU$ into $\CZ$, where,  for any control $\us\in\CU_R$, with 
associated state $(\phis,\mus,\ws )=:\CS(\us)$, the
first and second Fr\'echet derivatives $\CS'(\us)\in {\cal L}(\CU,\CZ)$ and $\CS''(\us)\in{\cal L}(\CU,{\cal L}(\CU,\CZ))$
are given as follows:

\vspace{1mm}\noindent
(i) \,\,\,\,For any  increment $h\in\CU$, $(\xi,\eta, v ):=\CS'(\us)[h]\in \CZ$ is the unique solution to the linearized problem
\begin{align}
\label{ls1}
&\dt \xi -\Delta \eta = 0 &&\mbox{a.e. in }\,Q,\\
\label{ls2}
&\tau \dt \xi -\Delta \xi - \eta - v = - f''(\phis)\xi &&\mbox{a.e. in }\,Q,\\
\label{ls3}
&\gamma \dt v + v = h  &&\mbox{a.e. in }\,Q,\\
\label{ls4}
&\dn \eta = \dn\xi = 0  &&\mbox{a.e. on }\,\Sigma,\\
\label{ls5}
&\xi(0)=0,\,\quad v(0)=0 &&\mbox{a.e. in }\,\Omega.
\end{align}
\Accorpa\Linear ls1 ls5

\vspace{1mm}
\noindent
(ii) \,\, For any pair of increments $h,k \in \CU$, $(\psi,\nu, z ):=\CS''(\us)[h,k]\in \CZ$ is the unique solution to the
bilinearized problem
\begin{align}
\label{bilin1}
&\dt \psi -\Delta \nu = 0 &&\mbox{a.e. in }\,Q,\\
\label{bilin2}
&\tau \dt \psi -\Delta \psi  - \nu - z = - f''(\phis)\psi - f'''(\phis)\xi^h\xi^k  &&\mbox{a.e. in }\,Q,\\
\label{bilin3}
&\gamma \dt z + z = 0  &&\mbox{a.e. in }\,Q,\\
\label{bilin4}
&\dn \nu = \dn\psi = 0 &&\mbox{a.e. on }\,\Sigma,\\
\label{bilin5}
&\psi(0)=0,\,\quad z(0)=0 && \mbox{a.e. in }\,\Omega,
\end{align}
\Accorpa\Bilinear bilin1 bilin5
where $\,(\xi^h,\eta^h, v^h ):=\CS'(\us)[h]\,$ and $\,(\xi^k,\eta^k, v^k ):=\CS'(\us)[k]$. We immediately note that the third component $z$ of the solution $(\psi,\nu, z )$ to \Bilinear\  fulfills
$z= 0$ a.e.\ in $Q$ due to \eqref{bilin3} and \eqref{bilin5}. 

\vspace*{3mm}
Now, fix some values $r_*,r^*$ such that
\Beq
\label{rstar}
-1<r_*<r_-<r_+<r^*<1,
\Eeq
with the constants $r_-,r_+$ introduced in \eqref{separation}. We then consider the set 
\begin{align}
\label{defPhi}
\Phi:=\biggl\{(\phi,\mu, w )\in\CZ:\  r_* < \infess_{(x,t)\in Q}\,\phi(x,t)
\leq \supess_{(x,t)\in Q}\,\phi(x,t) <r^* \biggr\},
\end{align}
which is obviously an open subset of $\CZ$. 
Notice that the functions in $\CZ$ are measurable and bounded, so that essential infimum and supremum are well defined.

We now prove an auxiliary result for the linear initial-boundary value problem
\begin{align}
\label{aux1}
&\dt \phi -\Delta \mu = 0 &&\mbox{a.e. in }\,Q,\\
\label{aux2}
&\tau \dt \phi -\Delta \phi - \mu - w = - \lambda_1 f''(\phis)\phi +\lambda_2 g &&\mbox{a.e. in }\,Q,\\
\label{aux3}
&\gamma \dt w + w = \lambda_3 h  &&\mbox{a.e. in }\,Q,\\
\label{aux4}
&\dn \mu = \dn\phi = 0 &&\mbox{a.e. on }\,\Sigma,\\
\label{aux5}
&\phi(0)=\lambda_4 \phi_0 ,\,\quad w(0)= \lambda_4 w_0 &&\mbox{a.e. in }\,\Omega, 
\end{align}
\Accorpa\Aux aux1 aux5
which for $\lambda_1=\lambda_3=1$ and $\lambda_2=\lambda_4=0$ coincides with the linearization \Linear\ of the state system. 
For convenience, we introduce the following Banach spaces for the initial data:
\begin{align}
\label{defN2}
&{\cal N}_2\,\,:=\{(\phi_0,w_0) : \ \phi_0\in W, \ w_0\in L^2(\Omega)\},\\
\label{defNinf}
&{\cal N}_\infty:=\{(\phi_0,w_0): \ \phi_0\in W, \ w_0\in L^\infty(\Omega)\},
\end{align}
equipped with their natural norms. We then have the following result.

\Blem
\label{Lem1}
Assume that $\lambda_1,\lambda_2,\lambda_3,\lambda_4\in\{0,1\}$ are given and that the assumptions {\bf (A1)}--{\bf (A3)} are fulfilled.
Moreover, let $(\us,(\phis,\mus,\ws ))\in \CU_R\times\Phi$ be arbitrary. Then the following holds true:\\[1mm]
{\rm (i)} \,\,\,\,The system \Aux\ has, for every   
$g, \, h\in \L2H$ and  $(\phi_0,w_0)\in{\cal N}_2$, a unique solution $(\phi,\mu,w) \in \CX$, and the linear mapping $\,(g,h, (\phi_0,w_0))
\mapsto (\phi,\mu,w)\,$ is continuous from $\,L^2(0,T;H) \times \L2H
\times{\cal N}_2\,$ into $\CX$.\\[1mm]
{\rm (ii)} \,\,If, in addition, $g\in H^1(0,T;H)$, then $(\phi,\mu,w)\in \CY$, and the mapping $\,(g,h, (\phi_0,w_0))
\linebreak \mapsto (\phi,\mu,w)\,$ is continuous from $\,H^1(0,T;H) \times \L2H
\times{\cal N}_2\,$ into $\CY$.\\[1mm]  
{\rm (iii)} If $\,g\in H^1(0,T;H) \cap \CU$, $h \in\CU$ and $(\phi_0,w_{0})\in {\cal N}_\infty$, 
then $(\phi,\mu,w) \in \CZ$, and the linear mapping 
$\,(g,h, (\phi_0,w_0))\mapsto (\phi,\mu,w)\,$
is continuous from $\,(\H1 H \cap \CU) \times \CU \times{\cal N}_\infty\,$ into~$\CZ$. 
\Elem

\Bdim
At first, arguing along the lines of the first part of the proof of Theorem~\ref{Teo1} for the state system, it is a standard 
matter to show that \Aux\ has a unique strong solution $(\phi,\mu,w) \in \CX$ for given data 
$g,\, h\in \L2H$, and
$(\phi_0,w_{0})\in {\cal N}_2$. Indeed, the existence can be proved via an appropriate Faedo--Galerkin 
approximation for which a priori estimates and a passage-to-the-limit process are performed. The uniqueness proof is 
actually simple due to the linearity of the problem. 
In order not to overload the exposition, we 
avoid writing the Faedo--Galerkin scheme here and just give the corresponding a priori estimates formally.  
To this end, we introduce the constants
\begin{align}
\label{defM1}
M_1&:=\lambda_2\,\|g\|_{\L2 H} +\lambda_3\,\|h\|_{\L2H}\,+\,\lambda_4\,\|(\phi_0,w_0)\|_{\CN_2},\\
\label{defM2}
M_2&:= \lambda_2\,\|g\|_{\H1 H}+ \lambda_3\,\|h\|_{\L2H}\,+\,\lambda_4\,\|(\phi_0,w_0)\|_{\CN_2},\\
\label{defM3}
M_3&:= \lambda_2\,\|g\|_{H^1(0,T;H)\cap \CU} +\lambda_3 \,\|h\|_{\CU}\,+\,\lambda_4\,\|(\phi_0,w_0)\|_{\CN_\infty}.
\end{align}

\begin{step}
Proof of (i)

Let $g\in L^2(0,T;H)$, $h\in \L2H$, and $(\phi_0,w_0)\in {\cal N}_2$ 
be given. We derive a series of a priori estimates for the corresponding solution $(\phi,\mu,w)$, where the 
constants $C>0$ are independent of the constants $M_1,M_2,M_3$.
 
We first add $\phi$ to both sides of \eqref{aux2} and then test \eqref{aux1} by $\mu$, \eqref{aux2} by $\dt \phi$, 
and \eqref{aux3} by $\dt w$. Addition of the resulting identities and integration over $(0,t) $ leads to a cancellation 
of terms, and it results that
\begin{align}
&\iQt |\nabla\mu|^2 +\tau \iQt |\dt \phi|^2 
+ \frac12 \| \phi(t)\|_V^2 + \iQt \gamma |\dt w|^2 +\frac12 \normaH{w(t)}^2  
\non \\
&=  - \lambda_1 \iQt  f''(\phis) \phi \, \dt\phi + \lambda_2 \iQt  g \, \dt\phi 
+ \lambda_3 \iQt  h\,  \dt w
\non \\
&\quad{} + \frac{\lambda_4^{\,2}}2 \| \phi_0\|_V^2  +
\frac{\lambda_4^{\,2}}2 \| w_0\|_H^2 + \iQt (\phi +w ) \dt \phi. 
\label{pier23}
\end{align}
Note that for the first term on the \rhs\ we can apply the bounds in Corollary~\ref{Cor1}, where now the constant $K_2$ depends on $r_*,\, r^*$ (in place of $r_-,\, r_+$). Then,
by the Young inequality we see that 
\begin{align*}
& - \lambda_1 \iQt  f''(\phis) \phi \, \dt\phi + \lambda_2 \iQt  g \, \dt\phi 
 + \iQt (\phi +w ) \dt \phi
\\
&\leq \frac \tau 2 \iQt |\dt \phi|^2 + C\bigl(\lambda_1^{\,2} +1 \bigr)\iQt
|\phi|^2  + C\lambda_2^{\,2} \iQt |g|^2  + C\iQt |w|^2 .
\end{align*} 
Similarly, we have that
\begin{align*}
 \lambda_3 \iQt  h\,  \dt w
\leq \frac{\gamma_0}2 \iQt |\dt w|^2 + C \lambda_3^{\, 2} \iQt |h|^2 .
\end{align*}
Then, using the fact that  $\iQt \gamma |\dt w|^2 \geq \gamma_0  \iQt |\dt w|^2$, and applying the above inequalities 
in \eqref{pier23}, we can infer from Gronwall's lemma that
\Beq
\label{pier24}
\norma{\nabla\mu}_{\revis{\L2 {H}^3}} + \norma{\phi}_{\H1H\cap \L\infty V} +\norma{w }_{\H1H} \leq C M_1.
\Eeq

Testing now \eqref{aux2} by $1/|\Omega| $, and integrating by parts with the help 
of \eqref{aux4}, we easily find, by comparison of terms and thanks to~\eqref{pier24}, that
\Beq
\label{pier25}
\norma{\overline\mu}_{L^2(0,T)}  \leq C M_1.
\Eeq

Next, in view of \eqref{pier24} and \eqref{pier25}, we can infer from the Poincar\'e--Wirtinger inequality 
that $\norma{\mu}_{\L2 V}
\leq C M_1$. Therefore, thanks to \eqref{aux1}, \eqref{aux4}, and elliptic regularity, it holds that 
$\norma{\Delta\mu}_{\L2 H}  \leq C M_1$ and 
\Beq
\label{pier26}
\norma{\mu}_{\L2 W}  \leq C M_1.
\Eeq
The same argument, this time applied to \eqref{aux2}, leads to  $\norma{\Delta\phi}_{\L2 H}  \leq C M_1$ and
\Beq
\label{pier27}
\norma{\phi}_{\L2 W}  \leq C M_1.
\Eeq

From the above estimates it follows that $(\phi,\mu, w)$ belongs to $\CX$ and, at the same time,  
that the continuity property asserted in (i) is valid. Assertion (i) is thus shown.
\end{step}

\begin{step}
Proof of (ii)

Assume now that $g\in \H1 H$. We then may differentiate \eqref{aux2} with respect to time and 
test by $\dt \phi$, then we add the resultant to \eqref{aux2} tested by $\dt \mu$. Again, we have a 
cancellation of terms, and by integration we obtain that
\begin{align}
&\frac12 \iO |\nabla\mu (t) |^2 + \frac\tau 2 \iO |\dt \phi(t) |^2 + \iQt |\nabla \dt \phi|^2 
\non
\\
&\leq  \frac12 \iO |\nabla\mu (0) |^2 + \frac\tau 2 \iO |\dt \phi(0) |^2
+ \iQt \dt w \, \dt \phi 
\non
\\
&\quad{}
 - \lambda_1 \iQt  f'''(\phis) \dt \phis \, \phi\, \dt\phi 
  - \lambda_1 \iQt  f''(\phis) |\dt \phi|^2 
 + \lambda_2 \iQt \dt g \, \dt\phi .
\label{pier28}
\end{align}
Concerning the first two terms on the \rhs, we can argue as in \eqref{pier11}--\eqref{pier13}, 
by reading \eqref{aux1} and \eqref{aux2} at the initial time $t=0$ and exploiting \eqref{aux5}. 
Since $\tau \dt \phi(0)$ is equal to $ \mu(0) $ plus the quantity 
$$
\lambda_4 (\Delta \phi_0  + w_0  - \lambda_1 f''(\phis(0) )\phi_0)  +\lambda_2 g(0),
$$
which is bounded in $H$ by $C M_2$ (see~\eqref{defPhi}, \eqref{defN2}, and note that 
$\normaH{g(0)} \leq C \norma{g}_{\H1 H}$ \,and $M_1\le M_2$), we infer that 
\begin{align*}
 \frac12 \iO |\nabla\mu (0) |^2 + \frac\tau 2 \iO |\dt \phi(0) |^2 \leq C M_2^2.
\end{align*}
Thanks to \eqref{pier24} and the bounds~\eqref{ssbound2} in Corollary~\ref{Cor1}, 
all of the other terms on the \rhs\ of \eqref{pier28} are easily under control except the fourth,
for which we argue as follows:  
\begin{align} 
&- \lambda_1 \iQt  f'''(\phis) \dt \phis \, \phi \, \dt\phi \leq 
C \lambda_1 \int_0^t \normaH{\dt \phis(s)} 
\norma{ \phi(s) }_{L^4(\Omega)} 
\norma{\dt \phi(s) }_{L^4(\Omega)}\,ds \non
\\
&\quad{} \leq C\lambda_1 \norma{\dt \phis}_{\L2 H} \norma{\phi}_{\L \infty V} 
\,\norma{\dt \phi(s) }_{L^2(0,t;V)} \non
\\
&\quad{} \leq \frac12 \iQt \bigl( |\dt \phi|^2 + |\nabla \dt \phi|^2 \bigr)  + C M_2^2,
\non
\end{align}
where we exploited the continuity of the embedding $V \subset L^4(\Omega)$. 
Then, combining the estimates above, we deduce from \eqref{pier28} and \eqref{pier24} that 
\Beq
\label{pier29}
\norma{\nabla\mu}_{\revis{\L\infty {H}^3}} + \norma{\dt \phi}_{\L{\infty}H\cap \L2 V}  \leq C M_2.
\Eeq
Now, we can repeat the comparison arguments used in \eqref{pier25}--\eqref{pier27} and conclude, in this order, that 
\begin{align}
\label{pier30}
&\norma{\overline\mu}_{L^\infty(0,T)}  \leq C M_2,
\\
\label{pier31}
&\norma{\mu}_{\L\infty W\cap \L2 {H^3(\Omega)}}  \leq C M_2,
\\
\label{pier32}
&\norma{\phi}_{\L\infty W}  \leq C M_2.
\end{align}
With these estimates, we have shown that $(\phi,\mu,w)\in\CY$ and that
$\|(\phi,\mu,w)\|_{\CY}\,\le\,C M_2,$
which concludes the proof of assertion (ii).  
\end{step}

\begin{step}
Proof of (iii)

Assume now that $g\in H^1(0,T;H)\cap\CU$, $h\in\CU$ and $(\phi_0,w_0)\in\CN_\infty$. Since $M_2\le CM_3$, we
then have from (ii) that $\|(\phi,\mu,w)\|_{\CY}\le CM_3$. Owing to the continuity of the embedding 
$W\subset L^\infty (\Omega)$ and to the fact that $w$ can be explicitly written as (cf.~\eqref{ssvar5})
\Beq
\non
w(x,t)= \lambda_4 w_0(x) \exp (-t/\gamma(x)) + \int_0^t \lambda_3 h(x,s)\exp(-(t-s)/\gamma(x) ) ds , \quad  (x,t) \in Q,
\Eeq
it is readily verified that $(\phi,\mu,w) $ belongs 
to $\CZ$, and, moreover, that  \,$\|(\phi,\mu,w)\|_{\CZ}\le C M_3$. 
This concludes the proof of the lemma.
\end{step}
\Edim

Now, having proved Lemma~\ref{Lem1}, we can prepare for the application of the implicit function
theorem. 
We consider two auxiliary linear initial-boundary value problems. The first is given by
\begin{align}
\label{oneaux1}
&\dt \phi -\Delta \mu = 0 &&\mbox{a.e. in }\,Q,\\
\label{oneaux2}
&\tau \dt \phi -\Delta \phi - \mu - w = g &&\mbox{a.e. in }\,Q,\\
\label{oneaux3}
&\gamma \dt w + w =  h  &&\mbox{a.e. in }\,Q,\\
\label{oneaux4}
&\dn \mu = \dn\phi = 0  &&\mbox{a.e. on }\,\Sigma,\\
\label{oneaux5}
&\phi(0)=0 ,\,\quad w(0)= 0 &&\mbox{a.e. in }\,\Omega,
\end{align}
\Accorpa\Oneaux oneaux1 oneaux5
and is obtained from \Aux\  for $\lambda_1=\lambda_4=0, \, \lambda_2=\lambda_3 = 1$.
Thanks to Lemma~\ref{Lem1}, the problem~\Oneaux\ admits for each $\,(g,h) \in (\H1 H \cap \CU) \times \CU \,$ a unique solution 
$ (\phi,\mu,w)\in \CZ $, and the associated linear mapping
$${\cal G}_1:(\H1 H \cap \CU) \times \CU \to \CZ,\quad  (g,h)\mapsto (\phi,\mu,w),$$ 
is continuous. The second system reads
\begin{align}
\label{twoaux1}
&\dt \phi -\Delta \mu = 0 &&\mbox{a.e. in }\,Q,\\
\label{twoaux2}
&\tau \dt \phi -\Delta \phi - \mu - w = 0 &&\mbox{a.e. in }\,Q,\\
\label{twoaux3}
&\gamma \dt w + w = 0  &&\mbox{a.e. in }\,Q,\\
\label{twoaux4}
&\dn \mu = \dn\phi = 0  &&\mbox{a.e. on }\,\Sigma,\\
\label{twoaux5}
&\phi(0)=\phi_0 ,\,\quad w(0)= w_0 &&\mbox{a.e. in }\,\Omega, 
\end{align}
\Accorpa\Twoaux twoaux1 twoaux5
and results from \Aux\ for $\lambda_1=\lambda_2=\lambda_3 =0,\lambda_4=1$.
For each $(\phi_0, w_0) \in \CN_\infty$, the problem~\Twoaux\ has a unique solution 
$ (\phi,\mu,w)\in \CZ $, and the associated mapping
$${\cal G}_2:\CN_\infty \to \CZ,\quad  (\phi_0, w_0)\mapsto (\phi,\mu,w),$$ 
is linear and continuous as well.
In addition, we define on the open set ${\cal A}:=({\cal U}_R\times\Phi)
\subset ({\cal U}\times\CZ)$ the nonlinear mapping
\begin{align}
\label{defG3}
&{\cal G}_3 :{\cal A}\to (\H1 H \cap \CU) \times \CU, \quad (u, (\phi,\mu,w)) \mapsto (-f'(\phi), u)
\end{align}
as a mapping from $\CU \times \CZ$ to $(\H1 H \cap \CU) \times \CU$.
The solution $(\phi,\mu,w)$ to the nonlinear state equation \State\ is the sum
of the solution to the system~\Twoaux\ and of the solution to the system
\Oneaux, where $(g,h)$ is given by the pair  $(-f'(\phi), u)$. 

All this means that the state $(\phi,\mu,w)$ associated with the control $
u $ is the unique solution to the nonlinear equation
\begin{equation} 
\label{nonlineq}
(\phi,\mu,w)= {\cal G}_2 (\phi_0, w_0) + {\cal G}_1  
 \big({\cal G}_3 (u, (\phi,\mu,w)) \bigr).
\end{equation}
Let us now define  the nonlinear mapping  $\,{\cal F}:{\cal A}\to \CZ$,
\begin{align}
\label{defF}
 {\cal F}(u, (\phi,\mu,w))\,:=\,{\cal G}_2 (\phi_0, w_0) + {\cal G}_1  
 \big({\cal G}_3 (u, (\phi,\mu,w)) \bigr) - (\phi,\mu,w).
\end{align} 
With $ {\cal F}$, the state equation can be shortly written as
\begin{equation} \label{nonlineq2}
{\cal F}(u, (\phi,\mu,w))=(0,0,0).
\end{equation}
This equation just means that $(\phi,\mu,w)$ is a solution to the state system \eqref{ss1}--\eqref{ss5} such that 
$(u, (\phi,\mu,w))\in{\cal A}$. From Theorem~\ref{Teo1} we 
know that such a solution exists for every $u\in{\cal U}_R$. A fortiori, any such solution automatically enjoys the separation property 
\eqref{separation} and is uniquely determined.  

We are going to apply the implicit function theorem to the equation  \eqref{nonlineq2}. To this aim,
we need the differentiability of the mappings entering \eqref{defF}. In particular, we have to
show that the mapping $\,{\cal G}_3\,$ is twice continuously Fr\'echet differentiable in ${\cal U}_R \times \Phi$ 
as a mapping from $\CU \times \CZ$ into $(H^1(0,T;H)\cap \CU)\times \CU$. To this end,
we first observe 
that, thanks to the differentiability properties of the
involved Nemytskii operators (see, e.g., \cite[Thm.~4.22, {p.~229}]{Fredibuch}), 
$\,{\cal G}_3\,$ is twice continuously Fr\'echet differentiable in ${\cal U}_R \times \Phi$ 
as a mapping from $\CU \times \CZ$ into $\CU\times \CU$,
and for the first partial derivatives at any point 
$\,(\us,(\phis,\mus,\ws))\in {\cal A}$, and for all $u\in{\cal U}$
and $(\phi,\mu,w)\in \CZ$, we have the identities
\begin{align}
\label{Freu}
&D_{u}{\cal G}_3 (\us,(\phis,\mus,\ws))[u]\,=\,
(0,u), \non \\
& D_{(\phi,\mu,w)}{\cal G}_3 (\us,(\phis,\mus,\ws))[(\phi,\mu,w)] = (-f''(\phis)\phi, 0). 
\end{align}  
It remains to show the differentiability properties of the mapping $(\phi,\mu,w)\mapsto -f'(\phi)$ on $\Phi$ as a mapping
from $\CZ$ into $H^1(0,T;H)$. Now let $(\phis,\mus,\ws)\in\Phi$ be fixed. In view of the explicit form of ${\cal G}_3$,
for the first derivative it obviously suffices to show that 
\begin{equation}
\label{prosecco1}
\frac{\|f'(\phis+\phi)-f'(\phis)-f''(\phis)\phi\|_{H^1(0,T;H)}}{\|(\phi,\mu,w)\|_{\CZ}}\,\to\,0
\quad\mbox{as }\,\|(\phi,\mu,w)\|_{\CZ}\to 0.
\end{equation}
To this end, let in the following $(\phi,\mu,w)\in\CZ$ be such that $(\phis+\phi,\mus+\mu,\ws+w)\in\Phi$. We then observe that
\begin{align}
\label{prosecco2}
&f'(\phis+\phi)-f'(\phis)-f''(\phis)\phi=\int_0^1 (1-\tau)\,f'''(\phis+\tau\phi)\,d\tau\,\phi^2=:A\,\phi^2, \\
\label{prosecco3}
&\mbox{with }\,\,|A|\le K, \quad |\dt A|\le K(|\dt\phis|+|\dt\phi|), \quad\mbox{a.e. in $Q$},  
\end{align}
where, here and in the following, $K>0$ denotes generic constants that are independent of the choice of $(\phi,\mu,w)$.
We thus have to estimate
\begin{align}
\label{prosecco4}
&\|A\,\phi^2\|_{H^1(0,T;H)}^2 \,=\,\iint_Q |A|^2\,|\phi|^4 +\iint_Q |\dt(A\,\phi^2)|^2 \, =: I_1+I_2.
\end{align}
Owing to \eqref{prosecco3}, we have
\begin{equation}
\label{prosecco5}
I_1\,\le\, K \,\|\phi\|_{\LiQ}^4\,\le\,K\,\|(\phi,\mu,w)\|_{\CZ}^4, 
\end{equation} 
as well as
\begin{align}
I_2\,&\le\,K\iint_Q(|\dt\phis|^2+|\dt\phi|^2)|\phi|^4\,+\,K\iint_Q |\phi|^2\,|\dt\phi|^2 \nonumber\\
&\le\,K\,\|\phi\|_{L^\infty(Q)}^4\,\left(1+\|\dt\phi\|^2_{L^2(0,T;H)}\right)
\,+\,K\,\|\phi\|_{L^\infty(Q)}^2\,\|\dt\phi\|^2_{L^2(0,T;H)} \nonumber\\
&\le K\,\|(\phi,\mu,w)\|_{\CZ}^4\left(1+\|(\phi,\mu,w)\|_{\CZ}^2\right). 
\end{align}
The validity of \eqref{prosecco1} is thus shown. 
The arguments for the second derivative and its continuity are quite similar, requiring only straightforward, albeit 
lengthy, calculations. To keep the paper at a reasonable length, we leave them to the interested reader. We just
remark at this place that the regularity requirement  $f_1, f_2\in C^5(-1,1)$ comes into play during the proof of the
continuity of the second derivative.

At this point, we introduce some abbreviating notation. We set
\begin{align*}
\by:=(\phi,\mu,w),\quad \bys:=(\phis,\mus,\ws), \quad \mathbf{0}=(0,0,0).
\end{align*}
Using the above differentiability results, we obtain from the chain rule  that ${\cal F}$ 
is twice continuously Fr\'echet differentiable in ${\cal U}_R\times\Phi$ as a mapping from 
$\CU\times {\cal Z}$	into $\CZ$, with the first-order partial derivatives
\begin{align}
D_u{\cal F}(\us,\bys)\,=\,{\cal G}_1\circ D_u{\cal 
 G}_3(\us,\bys), \quad D_{\by}{\cal F}(\us,\bys)\,=\,{\cal G}_1\circ
D_{\by}{\cal G}_3(\us,\bys)-I_{\CZ},
\label{DFy}
\end{align}
where $\,I_{\CZ}\,$ denotes the identity mapping on $\,\CZ$.

We want to prove the differentiability of the control-to-state mapping $ u \mapsto  {\bf y}$  defined implicitly 
by the equation $\,{\cal F}(u, {\bf y})=\mathbf{0}$, using the implicit function
theorem. Now let $u^*\in  {\cal U}_R$ be given and ${\bf y^*}={\cal S}(u^*)$. We need to show that the linear
and continuous operator $\,D_{{\bf y}}{\cal F}(u^*,{\bf y^*})$ is a topological isomorphism from $\CZ$ into itself.
 
To this end, let ${\bf v}=(v_1,v_2,v_3)\in\CZ$ be arbitrary. Then the identity $\,D_{\bf y}{\cal F}(u^*,{\bf y^*})[{\bf y}]={\bf v}\,$ 
just means that $\,{\cal G}_1\left(D_{\bf y}{\cal G}_3(u^*,{\bf y^*})[\bf y]\right)-{\bf y}={\bf v}$, which is equivalent to saying that   
\begin{equation*}
{\bf q}\,:=\,{\bf y}+{\bf v}= {\cal G}_1\left(D_{\bf y}{\cal G}_3(u^*,{\bf y^*})[{\bf q}]\right)
-{\cal G}_1\left(D_{\bf y}{\cal G}_3(u^*,{\bf y^*})[{\bf v}]\right).
\end{equation*}
The latter identity means that ${\bf q}$ is a solution to \eqref{aux1}--\eqref{aux3} for $\lambda_1=\lambda_2=\lambda_3=1,
\lambda_4=0$, with the specification $(g,h)=-D_{\bf y}{\cal G}_3(u^*,{\bf y^*})
[{\bf v}]=(f''(\phis)v_1,0)\in (H^1(0,T;H)\cap \CU)\times {\cal U}$. By Lemma~\ref{Lem1}, such a solution ${\bf q}\in\CZ$ exists 
and is uniquely 
determined, which shows that  $\,D_{\bf y}{\cal F}(u^*,{\bf y^*})$ is surjective. 
At the same time, taking ${\bf v}=\mathbf{0}$, we see that the equation 
$D_{{\bf y}}{\cal F}(u^*,{\bf y^*})[{\bf y}]=\mathbf{0}$ just
means that ${\bf y}$ is the unique solution to \eqref{aux1}--\eqref{aux3} for $\lambda_1=1,\,\,\lambda_2=\lambda_3=\lambda_4=0$. 
Obviously, ${\bf y}=\mathbf{0}$, which implies that $D_{{\bf y}}{\cal F}(u^*,{\bf y^*})$ is
also injective and thus, by the open mapping principle, a topological isomorphism from $\CZ$ into itself. 

We may therefore infer from the implicit function theorem (cf., e.g., \cite[Thms. 4.7.1 and 5.4.5]{Cartan} or \cite[10.2.1]{Dieu})     
that the control-to-state mapping $\CS$ is twice continuously Fr\'echet differentiable  in ${\cal U}_R$ 
as a mapping from $\CU$ into $\CZ$. More precisely, we obtain the following result.

\Bthm
\label{THM:FRECHET}
Suppose that the conditions {\bf (A1)}--{\bf (A3)} and $w_0\in L^\infty (\Omega)$ are fulfilled. 
Then the control-to-state operator
$\,\CS\,$ is twice continuously Fr\'echet differentiable in $\CU_R$ as a mapping from $\CU$ into $\CZ$.  
Moreover, for every $\us\in\CU_R$ and $h,k\in\CU$, the 
functions $(\xi,\eta,v)=\CS'(\us)[h]\in\CZ$ and $(\psi,\nu,z)=\CS''(\us)[h,k]\in\CZ$
are the unique solutions to the linearized system \Linear\ and the bilinearized system \Bilinear, respectively.
\Ethm
\Bdim
Let $\us\in\CU_R$ be arbitrary and $\bys=\CS(\us)$. The existence of $\CS'(\us)$ and $\CS''(\us)$, and their continuous dependence on $\us$, 
were shown above,
and differentiation of the identity ${\cal F}(u,\CS(u))=\mathbf{0}$ at $\us$ yields that
$$
D_{\by}{\cal F}(\us,\bys)\circ \CS'(\us)+ D_u{\cal F}(\us,\bys)=0.
$$ 
Now let $(\xi,\eta,v)=\CS'(\us)[h]$, where $h\in\CU$ is arbitrary. Then, by the above identity, and using
\eqref{DFy} and \eqref{Freu},
$$
(\xi,\eta,v)={\cal G}_1\bigl(D_{\by}{\cal G}_3(\us,\bys)[(\xi,\eta,v)]\,+\,D_u{\cal G}_3(\us,\bys)[h]\bigr)
={\cal G}_1((-f''(\phis)\xi,h)),
$$
and it easily follows from the definition of ${\cal G}_1$ that $(\xi,\eta,v)$ indeed coincides
with the unique solution to \Linear\ which, by Lemma~\ref{Lem1},(iii), belongs to $\CZ$. 

The calculation of the form of the second derivative $\CS''(\us)$ is not given here in order to keep the 
exposition at a reasonable length. We just mention that the arguments employed in \cite[Sect.~5.7]{Fredibuch}
for  a semilinear heat conduction problem carry over to our situation with only minor changes, leading to the
conclusion that $(\psi,\nu,z)=\CS''(\us)[h,k]$ indeed solves the system \Bilinear. Now observe that the system
\Bilinear\ is of the form \eqref{aux1}--\eqref{aux5} with $\lambda_1=\lambda_2=1, \,\,\lambda_3=\lambda_4=0$,
and $g:=-f'''(\phis)\xi^h\xi^k$. It is not difficult to show that $g\in \H1 H \cap \CU$, and Lemma~\ref{Lem1},(iii)
yields that $(\psi,\nu,z)\in \CZ$. 
\Edim

\Brem
It is worth noting that for the argumentation used above the actual value of the constant $R>0$ defining $\CU_R$ did not matter.
It therefore follows that $\CS$ is twice continuously Fr\'echet differentiable as a mapping from $\CU$ to $\CZ$ on the entire
space $\CU$.
\Erem
\Brem
\label{Spextended}
In view of the continuous embedding $\CZ\subset\CY$, the control-to-state mapping $\CS$ is also Fr\'echet differentiable
from $\CU$ to $\CY$ with the same expression for the Fr\'echet derivative, now regarded as an element of ${\cal L}(\CU,\CY)$. 
As $\CU$ is dense in $\L2H$, the operator $\CS'(\us)\in {\cal L}(\CU,\CY)$ can be extended in the standard way to an operator 
belonging to ${\cal L}(\L2H,\CY)$. We still denote the extended operator by $\CS'(\us)$, where we underline that it coincides 
with a Fr\'echet derivative only on $\CU$ and not on $\L2H$. However, it is readily seen by a density argument that $(\xi,\eta,v)
=\CS'(\us)[h]$ coincides also for $h\in \L2 H$ with the solution to \Linear. Analogously, the second Fr\'echet derivative $\CS''(\us)$ can
be continuously extended, which leads to an element of the space ${\cal L}(L^2(0,T;H),{\cal L}(L^2(0,T;H),\CY))$ that is 
still denoted by $\CS''(\us)$. Again, $(\psi,\nu,z)=\CS''(\us)$ $[h,k]$ solves \Bilinear\ also for $h,k\in\L2H$.
For the extensions, we have the following result.  
\Bcor
Let {\bf (A1)--(A3)} and $w_0\in L^\infty (\Omega)$ be fulfilled, and let $\us\in \CU_R$ be fixed. Then we have for every $h,k\in\L2H$ the estimates
\begin{equation}
\label{extension}
\|\CS'(\us)[h]\|_{\CY}\, \le K_6\,\|h\|_{L^2(0,T;H)}, \quad \|\CS''(\us)[h,k]\|_{\CY}\,\le\,K_6\,
\|h\|_{L^2(0,T;H)}\,\|k\|_{L^2(0,T;H)},
\end{equation}
with a constant $K_6>0$ that depends only on $R$ and the data. 
\Ecor
\Bdim
First note that $(\xi,\eta,v)=\CS'(\us)[h]$ solves the system \Linear, which is of the form \Aux\ with $\lambda_1=\lambda_2
=\lambda_3=1, \,\,\lambda_4=0$ and $g=0$. Therefore the first inequality in \eqref{extension}
follows directly from Lemma~\ref{Lem1},(ii). Next, we 
observe that the system \Bilinear\ is also of the form \Aux, but this time with $\lambda_1=\lambda_2=1,\,\,\lambda_3=\lambda_4=0$,
and $g=-f'''(\phis)\xi^h\xi^k$. Hence, also the second inequality in \eqref{extension} will follow from Lemma~\ref{Lem1},(ii)
once we can show that  
\begin{equation}
\label{Jojo}
\|g\|_{H^1(0,T;H)}\le \widehat C\,\|h\|_{L^2(0,T;H)}\,\|k\|_{L^2(0,T;H)}
\end{equation}
with some $\widehat C>0$ that only depends on $R$ and the data. Now recall the definition \eqref{defY} of $\CY$ and the fact that the
first estimate in \eqref{extension} is already shown. We thus have 
\begin{align*}   
\|g\|_{L^2(0,T;H)}^2 \,&\le \,C\,\|\xi^h\|^2_{\LiQ}\|\xi^k\|^2_{\LiQ}\,\le\,C\|\CS'(\us)[h]\|^2_{\CY}\,\|\CS'(\us)[k]\|_{\CY}^2\\
&\le\,C\,\|h\|^2_{L^2(0,T;H)}\,\|k\|^2_{L^2(0,T;H)}\,.
\end{align*}
Moreover,
\begin{align*}
&\|\dt g\|_{L^2(0,T;H)}^2\,\le\,C\iint_Q |\dt\phis\,\xi^h\,\xi^k|^2\,+\,C\iint_Q \bigl(|\dt\xi^h|^2 |\xi^k|^2\,+
\,|\xi^h|^2 |\dt\xi^k|^2\bigr)\\		
&\le\,C\,\|\xi^h\|^2_{\LiQ}\,\|\xi^k\|^2_{\LiQ} \,+\,C\bigl(\|\dt\xi^h\|^2_{L^2(Q)}\,\|\xi^k\|^2_{\LiQ}
\,+\,\|\xi^h\|^2_{\LiQ}\,\|\dt\xi^k\|^2_{L^2(Q)}\bigr)\\
&\le\,C\|\CS'(\us)[h]\|^2_{\CY}\,\|\CS'(\us)[k]\|_{\CY}^2\,\,\le\,\,C\,\|h\|^2_{L^2(0,T;H)}\,\|k\|^2_{L^2(0,T;H)}\,,
\end{align*}
which concludes the proof.
\Edim
\Erem

Next, we show a Lipschitz continuity property of the extensions of the derivatives that will prove crucial for
the derivation of second-order sufficient optimality conditions below. 

\Bthm 
The  mappings $\,\CU\to {\cal L}(L^2(0,T;H),\CY)$, $\,u \mapsto \CS'(u)$, and\, 
$\CU\to {\cal L}(L^2(0,T;\linebreak H),{\cal L}(\L2 H,\CY))$, $u \mapsto \CS''(u)$,\,  are 
Lipschitz continuous in the following sense: there exists a constant
$K_6>0$, which depends only on $R$ and the data, such that, for all controls $u_1,u_2 \in\CU_R$ and all 
increments $h,k \in\L2 H$, it holds that 
\begin{align}
\label{lip1} 
&\|(\CS'(u_1)-\CS'(u_2))[h]\|_{\CX}
\le\,K_6\,\|u_1-u_2\|_{\L2H}\,\|h\|_{\L2H}\,,\\[2mm]
\label{lip2}
&\|\left( \CS''  (u_1)- \CS'' (u_2)\right)[h,k]\|_{\CX}  
\le\,K_6\,\|u_1-u_2\|_{\L2H}\,\|h\|_{\L2H}\,\|k\|_{\L2H}\,.
\end{align}
\Ethm
\Bdim
Let $u_1,u_2\in\CU_R$ and $h,k\in\L2H$ be given. We put
\begin{align*}
&(\phi_i,\mu_i,w_i)=\CS(u_i), \quad (\xi_i^h,\eta_i^h,v_i^h)=\CS'(u_i)[h], \quad (\xi_i^k,\eta_i^k,v_i^k)=\CS'(u_i)[k],\\
&(\psi_i,\nu_i,z_i)=\CS''(u_i)[h,k], \quad\mbox{for \,$i=1,2$, as well as} \\
&(\xi^h,\eta^h,v^h)=(\xi_1^h-\xi_2^h,\eta_1^h-\eta_2^h,v_1^h-v_2^h), \quad 
(\xi^k,\eta^k,v^k)=(\xi_1^k-\xi_2^k,\eta_1^k-\eta_2^k,v_1^k-v_2^k).
\end{align*}
Then it is easily verified that the triple $(\xi^h,\eta^h,v^h)$ solves a system of the form \Linear, only that $h=0$ in this case and
that the \rhs\ of 
\eqref{ls2} is here replaced by the expression 
$$
\widetilde g:= -f''(\phi_1)\xi^h-(f''(\phi_1)-f''(\varphi_2))\xi_2^h\,.
$$
Now observe that this system is of the form \Aux\ with $\lambda_1=\lambda_2=1, \,\,\lambda_3=\lambda_4=0$, \,$\phis=\phi_1$, and
$\,g=-(f''(\phi_1)-f''(\varphi_2))\xi_2^h$. Therefore it follows from Lemma~\ref{Lem1},(i) that the inequality \eqref{lip1} is
valid provided we can show that 
\begin{equation}
\label{Jojo1}
\|g\|_{L^2(0,T;H)}\,\le\,C\,\|u_1-u_2\|_{L^2(0,T;H)}\,\|h\|_{L^2(0,T;H)}.
\end{equation}

Now, by \eqref{ssbound2}, it holds $\,|g| \le C|\phi_1-\phi_2||\xi_2^h|\,$ a.e. in $Q$. Invoking \eqref{contdep2} and 
\eqref{extension}, we therefore conclude that
\begin{align*}
\|g\|^2_{L^2(Q)}&\le\,C\iint_Q|\phi_1-\phi_2|^2\,|\xi_2^h|^2\,\le\,C\,\|\phi_1-\phi_2\|_{L^2(Q)}^2\,\|\xi_2^h\|^2_{\LiQ}\\
&\le \,C\|\CS(u_1)-\CS(u_2)\|_{\CX}^2\,\|\CS'(u_2)[h]\|^2_{\CY}\,\le\,C\,\|u_1-u_2\|_{L^2(Q)}^2\,\|h\|^2_{L^2(Q)}.
\end{align*}  
The inequality \eqref{lip1} is thus proved. To show the validity of \eqref{lip2}, we observe that the triple 
$$
(\psi,\nu,z)=(\psi_1-\psi_2,\nu_1-\nu_2,z_1-z_2)
$$
satisfies a system of the form \Bilinear, where this time the \rhs\ of \eqref{bilin2} is given by
$$
\widetilde g = - f''(\phi_1)\psi -(f''(\phi_1)-f''(\phi_2))\psi_2- (f'''(\phi_1)\xi_1^h\xi_1^k-f'''(\phi_2)\xi_2^h\xi_2^k).
$$
The system for $(\psi,\nu,z)$ is again of the form \Aux, this time with $\lambda_1=\lambda_2=1,\,\,\lambda_3=\lambda_4=0$, $\phis=\phi_1$, and 
$$g=-(f''(\phi_1)-f''(\phi_2))\psi_2- (f'''(\phi_1)\xi_1^h\xi_1^k-f'''(\phi_2)\xi_2^h\xi_2^k),$$
and, thanks to Lemma~\ref{Lem1},(i),  it suffices to show that 
\begin{equation}
\label{Jojo2}
\|g\|_{L^2(0,T;H)}\,\le\,C\,\|u_1-u_2\|_{L^2(0,T;H)}\,\|h\|_{L^2(0,T;H)}\,\|k\|_{L^2(0,T;H)}\,.
\end{equation}
Now observe that \eqref{ssbound2} yields that, almost everywhere in $Q$,
\begin{align*} 
|g|&\le \,C\bigl(|\phi_1-\phi_2||\psi_2|\,+\,|\phi_1-\phi_2||\xi_1^h||\xi_1^k|\,+\,|\xi_1^h-\xi_2^h||\xi_1^k|
\,+\,|\xi_2^h||\xi_1^k-\xi_2^k|\bigr).
\end{align*}          
Hence, by virtue of \eqref{contdep2}, \eqref{extension} and the already shown estimate \eqref{lip1},
\begin{align*}
\|g\|_{L^2(Q)}&\le\,C\bigl(\|\phi_1-\phi_2\|_{\L2H}\,\|\psi_2\|_{\LiQ}\,+\,\|\phi_1-\phi_2\|_{\L2H}\,\|\xi_1^h\|_{\LiQ}
\,\|\xi_1^k\|_{\LiQ}\\
&\qquad\quad +\,\|\xi^h\|_{\L2H}\,\|\xi_1^k\|_{\LiQ}\,+\,\|\xi_2^h\|_{\LiQ}\,\|\xi^k\|_{\L2H}\bigr)\\
&\le\,C\Big(\|\CS(u_1)-\CS(u_2)\|_{\CX}\,\bigl(\|\CS''(u_2)[h,k]\|_{\CY}\,+\,\|\CS'(u_1)[h]\|_{\CY}\,\|\CS'(u_1)[k]\|_{\CY}\bigr)\\
&\qquad\quad +\,\|(\CS'(u_1)-\CS'(u_2))[h]\|_{\CX}\,\|\CS'(u_1)[k]\|_{\CY} \\
&\qquad\quad +\,\|\CS'(u_2)[h]\|_{\CY}\,
\|(\CS'(u_1)-\CS'(u_2))[k]\|_{\CX}\Big)\\
&\le\,C\,\|u_1-u_2\|_{L^2(0,T;H)}\,\|h\|_{L^2(0,T;H)}\,\|k\|_{L^2(0,T;H)},
\end{align*}
which concludes the proof of the assertion.
\Edim

\section{The optimal control problem}
 
\setcounter{equation}{0}

In this section, we study the optimal control problem {\bf (CP)} with the cost functional \eqref{cost}. Besides the general
conditions {\bf (A1)}--{\bf (A3)} and $w_0\in L^\infty (\Omega)$, we  make the following 
assumptions:
\begin{description}
\item[(A4)] \,\,\juerg{It holds $b_1\ge 0$, $b_2\ge 0$, $b_3>0$,} and $\,\kappa>0$. 
\item[(A5)] \,\,The target functions satisfy $\phi_Q\in L^2(Q)$ and $\phi_\Omega \in V.$
\end{description}
We assume $\kappa > 0$ to include the effects of sparsity. By an obvious modification, the theory of second-order 
conditions remains valid for $\kappa = 0$.
\Brem
The assumption $\phi_\Omega \in V$ is  useful in order to have more regular solutions to the associated
adjoint system (see below). It is not overly restrictive in view of the continuous embedding \,$(H^1(0,T;H)\cap L^2(0,T;W))
\subset C^0([0,T];V)$\, which implies that $\phi(T)\in V$. 
\Erem

The following existence result can be shown with a standard argument that needs no repetition here 
(see, e.g., a similar result with proof in \cite[Thm.~4.1]{CGRS}). 

\Bthm
Suppose that {\bf (A1)}--{\bf (A5)} are fulfilled, and suppose that $\juerg{G}:L^2(Q)\to\erre$ is 
\juerg{nonnegative,} convex and continuous. Then the optimal control problem {\bf (CP)} admits a solution 
$\us\in\Uad$.
\Ethm

\subsection{The adjoint system}

In the following, we often denote by $\us\in\Uad$ a locally optimal control for {\bf (CP)} and by $(\phis,\mus,\ws)=\CS(\us)$
the associated state. Recall that a control $\us\in\Uad$ is called locally optimal in the sense of $L^p(Q)$ for some
$p\in[1,+\infty]$ if and only if there is some $\varepsilon>0$ such that \revis{$\,{\cal J}(\us,\CS_1(\us))\le {\cal J}(u,\CS_1(u))\,$} 
for all $u\in\Uad$ with $\|u-\us\|_{L^p(Q)}\le\varepsilon$. As can easily be seen, any locally optimal control in the sense of $L^p(Q)$ for \revis{some}
$1\le p<+\infty$ is also locally optimal in the sense of $L^\infty(Q)$. 

 The corresponding adjoint state system is formally given by:
\begin{align}
\label{adj1}
&-\dt( p+\tau q)-\Delta q + f''(\phis )q =b_1(\phis -\phi_Q) &&\mbox{a.e. in }\,Q,\\
\label{adj2}
&-\Delta p-q=0 &&\mbox{a.e. in }\,Q,\\
\label{adj3}
&-\gamma\dt r+r-q=0 &&\mbox{a.e. in }\,Q,\\
\label{adj4}
&\dn p=\dn q=0 &&\mbox{a.e. on }\,\Sigma,\\
\label{adj5}
&(p+\tau q)(T)=b_2(\phis(T)-\phi_\Omega), \quad r(T)=0 &&\mbox{a.e. in }\,\Omega.
\end{align}
\Accorpa\Adjoint adj1  adj5
We immediately observe that the system is decoupled in the sense that $\,r\,$ can be directly recovered from \eqref{adj3} with the
terminal condition $r(T)=0$ once $q$ is determined. Note also that the variational form of \eqref{adj1}, \eqref{adj2}, \eqref{adj4} 
is given by
\begin{align}
\label{wadj1}
& -\iO\dt(p+\tau q)\rho +\iO\nabla q\cdot\nabla\rho +\iO f''(\phis)q\rho
=b_1\iO (\phis-\phi_Q)\rho \nonumber\\
&\qquad \mbox{for a.e. $\,t\in(0,T)$\, and every }\,\rho\in V,\\
\label{wadj2}
&\iO \nabla p\cdot\nabla\rho=\iO q\rho \quad\mbox{for a.e. $\,t\in(0,T)$\, and every }\,\rho\in V. 
\end{align}
We have the following result. 
\Bthm
\label{Teo3}
Suppose that {\bf (A1)}--{\bf (A5)}, \revis{{}$\gamma \in W^{2,\infty}(\Omega)$ and{}} $w_0\in L^\infty (\Omega)$ are fulfilled, and let $\us \in \CU_R$ be a control with associated state $(\phis,\mus,\ws)$. 
Then the associated adjoint state system
has a unique strong solution $(\ps,\qs,\rs)$ with the regularity
\begin{align}
\label{regpq}
&\ps+\tau \qs\in H^1(0,T;H)\cap C^0([0,T];V)\cap L^2(0,T;W),\\
\label{regp}
&\ps\in L^2(0,T;W\cap H^4(\Omega)),\\
\label{regq}
&\qs\in L^2(0,T;W),\\
\label{regr}
&\rs\in H^1(0,T;W).
\end{align}
Moreover, there is a constant $K_7>0$, which depends only on $\,R\,$ and the data, such that 
\begin{align}
&\|\ps+\tau\qs\|_{H^1(0,T;H)\cap C^0([0,T];V)\cap L^2(0,T;H^2(\Omega))}\,+\,\|\ps\|_{L^2(0,T;H^4(\Omega))} \,+\,\|\qs\|_{L^2(0,T;H^2(\Omega))}
\nonumber\\
&\quad +\,\|\rs\|_{H^1(0,T;H^2(\Omega))}\,\le\, K_7\left(\|\phis-\phi_Q\|_{L^2(Q)} + \|\phis(T)-\phi_\Omega\|_V \right).\label{adj6}
\end{align}
\Ethm
\Bdim We solve the initial-boundary value problem given by \eqref{wadj1} and \eqref{wadj2} together with the first terminal condition
in \eqref{adj5} via a Faedo--Galerkin approximation. To this end, let $\{\lambda_j\}_{j\in\enne}$ and $\{e_j\}_{j\in\enne}$ denote 
the countable sets  of eigenvalues and eigenfunctions to the
elliptic eigenvalue problem $-\Delta e_j=\lambda_j e_j$ in $\Omega$, $\dn e_j=0$ on $\Gamma$, where the eigenfunctions are normalized by 
$\|e_j\|_{L^2(\Omega)}=1$ for $j\in\enne$. Then 
\begin{align*}
&0=\lambda_1<\lambda_2\le \lambda_3\le\ldots, \quad\lim_{j\to\infty}\lambda_j=+\infty,\nonumber\\
&\iO e_j e_k=\iO\nabla e_j\cdot\nabla e_k=0 \quad\mbox{for }\,j\not=k.
\end{align*} 
We now introduce the $n$-dimensional spaces $\,V_n:=\mbox{span} \{e_1,\ldots,e_n\}$, $n\in\enne$,
 and recall the well-known fact that $\bigcup_{n\in\enne} V_n$ is dense
in both $H$ and $V$.  
We make the ansatz
\begin{equation*}
p_n(x,t)=\sum_{j=1}^n p_{n,j}(t)e_j(x), \quad q_n(x,t)=\sum_{j=1}^n q_{n,j}(t)e_j(x), \quad\mbox{for }\,(x,t)\in\overline Q,
\end{equation*}
and look for functions $p_{n,j},\,\,q_{n,j}$ such that the identities \eqref{wadj1} and \eqref{wadj2} are fulfilled, where
$p,\,q$ are replaced by $p_n,\,q_n$ and the test functions $\rho$ are required to belong to $V_n$; moreover, we postulate 
the terminal condition $\,(p_n+\tau q_n)(T)= \Pi_n\bigl(b_2(\phis(T)-\phi_\Omega)\bigr)$, where $\Pi_n$ denotes the $H$-orthogonal projection
operator onto $V_n$. Observe that 
\begin{align}
\label{pin1}
&\Pi_n v=\sum_{j=1}^n (v,e_j)_H\, e_j \quad\mbox{and}\quad \|\Pi_n v\|_H\le\|v\|_H \quad\mbox{for all }\,v\in H,\\
\label{pin2}
&\|\nabla(\Pi_n v)\|_H \le\,\|v\|_V\quad\mbox{for all }\,v\in V.
\end{align}
Next, we choose $\rho=e_i$ in \eqref{wadj1} and \eqref{wadj2}, which leads to the system
\begin{align}
\label{ODE1}
&-\dt(p_{n,i}+\tau q_{n,i})+\lambda_i q_{n,i}+\iO f''(\phis) \sum_{j=1}^n q_{n,j}e_j e_i =b_1\iO (\phis-\phi_Q)e_i,
\nonumber\\
&\qquad \mbox{a.e. in }\,(0,T), \quad\mbox{for }\,1\le i\le n,\\
\label{ODE2}
&\lambda_i\,p_{n,i}=q_{n,i} \quad\mbox{a.e. in }\,(0,T), \quad\mbox{for }\,1\le i\le n.
\end{align}
In addition, testing the terminal condition by $e_i$, we find that
\begin{equation}
\label{ODE3}
(p_{n,i}+\tau q_{n,i})(T)= b_2\iO (\phis(T)-\phi_\Omega)e_i\quad\mbox{for }\,1\le i\le n.
\end{equation}
Now we substitute for $q_{n,i}$ from \eqref{ODE2} in \eqref{ODE1} and \eqref{ODE3}, which leads to an explicit backward Cauchy problem
for a nonhomogeneous linear ODE system in the unknowns $p_{n,1},\ldots,p_{n,n}$ whose coefficient functions and \rhs s all
belong to $L^2(0,T)$. Owing to Carath\'eodory's theorem, there exists a unique solution $(p_{n,1},\ldots,p_{n,n})\in H^1(0,T;\erre^n)$,
which, in turn, uniquely determines the solutions $p_n,q_n\in H^1(0,T;W)$ to the $n$-dimensio\-nal version of the
variational system \eqref{wadj1}, \eqref{wadj2}, together with the terminal condition 
$\,(p_n+\tau q_n)(T)= \Pi_n\bigl(b_2(\phis(T)-\phi_\Omega)\bigr)$. We now derive a number of a priori estimates for $p_n$ 
and $q_n$, where in the
following $C>0$ denotes constants that may depend on the data, but not on $n\in\enne$. When saying that we ``insert functions
in \eqref{wadj1} or \eqref{wadj2}'', we will always mean the $n$-dimensional versions of these variational equalities which are 
solved by~$p_n$~and~$q_n$. 

Let $n\in\enne$ now be fixed, and let 
\begin{equation}
\label{Mneu}
M:=\|\phis-\phi_Q\|_{L^2(Q)} + \|\phis(T)-\phi_\Omega\|_V.
\end{equation}
First, we insert $\rho=p_n+\tau q_n$ in \eqref{wadj1} and $\rho=q_n$ in \eqref{wadj2}, and subtract the resultants, noting that a 
cancellation of two terms occurs. Then we integrate over $(t,T)$, where $t\in[0,T)$ is arbitrary. Introducing the notation
$Q^t:=\Omega\times(t,T)$ for $t\in[0,T)$, and using Young's inequality, \eqref{ssbound2}, and the fact that 
$\|\Pi_n(b_2(\phis(T)-\phi_\Omega))\|_H
\,\le\,\|b_2(\phis(T)-\phi_\Omega)\|_H$ by \eqref{pin1}, we then obtain that
\begin{align}
&\frac 12\|(p_n+\tau q_n)(t)\|_H^2\,+\,\iint_{Q^t}|q_n|^2\,+\,\tau\iint_{Q^t}|\nabla q_n|^2\nonumber\\
&=\,\frac 12 \|\Pi_n(b_2(\phis(T)-\phi_\Omega))\|_H^2\,-\,\iint_{Q^t} f''(\phis)q_n(p_n+\tau q_n)\nonumber\\
&\qquad + \,b_1\iint_{Q^t} (\phis-\phi_Q)(p_n+\tau q_n)\nonumber\\
&\le \,CM^2\,+\,\frac 12\iint_{Q^t}|q_n|^2\,+\,C\iint_{Q^t}|p_n+\tau q_n|^2,\nonumber
\end{align}
and Gronwall's lemma yields that
\begin{equation}
\label{adjesti1}
\|p_n+\tau q_n\|_{L^\infty(0,T;H)}\,+\,\|q_n\|_{L^2(0,T;V)}\,\le\,CM.
\end{equation}
In addition, we conclude from \eqref{ODE2} that $\,-\Delta p_n=q_n\,$ a.e.~in $(0,T)$; since also $\,\dn p_n=0\,$ on $\Gamma$
a.e.~in $(0,T)$, we can therefore infer from \eqref{adjesti1} and elliptic regularity theory that $p_n\in L^2(0,T;W\cap H^3(\Omega))$
and also
\begin{equation}
\label{adjesti2}
\|p_n\|_{L^2(0,T;H^3(\Omega))}\,\le\,CM.
\end{equation}
Moreover, it readily follows from \eqref{adjesti1} and \eqref{adjesti2}, by comparison in \eqref{wadj1}, that
$$
\Big|\iint_Q -\dt(p_n+\tau q_n)\rho \Big|\,\le\,CM\|\rho\|_{L^2(0,T;V)} \quad\mbox{for all }\,\rho\in L^2(0,T;V_n).
$$  
A standard argument then yields that $\dt(p_n+\tau q_n)\in L^2(0,T;V^*)$ and
\begin{equation}
\label{adjesti3}
\|p_n+\tau q_n\|_{H^1(0,T;V^*)}\,\le\,CM.
\end{equation}

At this point, we can apply well-known weak and weak-star compactness arguments to 
conclude that there are functions $p^*,q^*$ such that as $n\to\infty$ 
(at first only for a suitable subsequence, but due to the uniqueness of the limit point eventually for the entire sequence) it holds
\begin{align}
\label{conpq}
p_n+\tau q_n&\to p^*+\tau q^*&&\mbox{weakly star in }\,H^1(0,T;V^*)\cap L^\infty(0,T;H),\\
\label{conp}
p_n &\to p^* &&\mbox{weakly in }\,L^2(0,T;H^3(\Omega)),\\
\label{conq}
q_n &\to q^* &&\mbox{weakly in }\,\L2 V.
\end{align} 
Besides, standard arguments (which need no repetition here) imply that the pair $(p^*,q^*)$ satisfies \eqref{wadj1}, \eqref{wadj2}, and
the terminal condition $(p^*+\tau q^*)(T)=b_2(\phis(T)-\phi_\Omega)$. Also, by virtue of the semicontinuity properties of norms, we infer that the estimates 
\eqref{adjesti1}--\eqref{adjesti3} are valid with $p_n,q_n$ replaced by $p^*,q^*$.
Moreover, from the linearity of \eqref{wadj1} and \eqref{wadj2} it readily
follows that the solution is unique.

As the next step, we now recover further regularity properties for $\ps,\qs$. To this end, we 
insert $\rho=-\Delta(p_n+\tau q_n)\in V_n$ in \eqref{wadj1} for fixed $n\in\enne$ and integrate with respect to time to obtain that 
\begin{align} 
\label{adjesti4}
&\frac 12 \iO |\nabla(p_n+\tau q_n)(t)|^2\,+\,\tau\iint_{Q^t}|\Delta q_n|^2\,=\,\frac 12\iO|\nabla (\Pi_n(b_2(\phis(T)-\phi_\Omega)))|^2
\nonumber\\
&\quad-\iint_{Q^t}\Delta q_n\Delta p_n\,-\iint_{Q^t}\bigl(b_1(\phis-\phi_Q)-f''(\phis)q_n\bigr)\Delta (p_n+\tau q_n).
\end{align}
Now recall that $v:=\phis(T)-\phi_\Omega\in V$, so that \eqref{pin2} can be applied. Therefore, applying Young's inequality
and the bounds \eqref{adjesti1} and \eqref{adjesti2} to the last two terms on the \rhs\ of \eqref{adjesti4}, we conclude that
the \revis{\lhs}\ is bounded by \revis{$\,CM^2$}. Therefore, by virtue of elliptic regularity,
\begin{equation}
\label{adjesti5}
\|p_n+\tau q_n\|_{L^\infty(0,T;V)}\,+\,\|q_n\|_{L^2(0,T;\Hdue)}\,\le\,CM,
\end{equation}
whence, since $-\Delta p_n=q_n$, and using elliptic regularity theory once more, 
\begin{equation}
\label{adjesti6}
\|p_n\|_{L^2(0,T;H^4(\Omega))}\,\le\,CM.
\end{equation}
Consequently, the limit points satisfy \eqref{regp} and \eqref{regq}. Moreover, we have $p^*+\tau q^*\in 
L^\infty(0,T;V)\cap L^2(0,T;\Hdue)$ and $q^*\in L^2(0,T;\Hdue)$ with the corresponding norm estimates. 
In addition, comparison in \eqref{wadj1}, using the already shown bounds,  yields that also $p^*+\tau q^*\in H^1(0,T;H)$ together
with the estimate
\begin{equation}
\label{adjesti7}
\|p^*+\tau q^*\|_{H^1(0,T;H)}\,\le\,CM.
\end{equation}
By virtue of the continuous embedding $\bigl(H^1(0,T;H)\cap L^2(0,T;\Hdue)\bigr)\subset C^0([0,T];V)$, then also
\begin{equation}
\label{adjesti8}
\|p^*+\tau q^*\|_{C^0([0,T];V)}\,\le\, CM.
\end{equation} 

It thus remains to show that $\,\,\|r^*\|_{H^1(0,T;\Hdue)}\le CM$ for the uniquely determined function $r^*$
satisfying \eqref{adj3} with $q=\qs$ and $r^*(T)=0$. \revis{Note that $r^*$ and each partial derivative of first and second order of $r^*$  solves a linear ODE with coefficients in $L^\infty (\Omega) $ and known terms in $\L2H$, with the main coefficient $\gamma$ being bounded from below by $\gamma_0 >0$, as from {\bf (A2)}.  Then the searched estimate is a} consequence of the estimate for 
$q^*$ following from \eqref{adjesti5}. This concludes the proof of the assertion.
\Edim 

The following continuous dependence result will be needed below during the proof of second-order sufficient optimality conditions.
\Bcor
\label{Cor4.4}
Suppose that {\bf (A1)}--{\bf (A5)}, \revis{$\gamma \in W^{2,\infty}(\Omega)$}  and $w_0\in L^\infty (\Omega)$ are fulfilled, 
and let, for $i=1,2$, $u_i\in\CU_R$ be given with the associated states
$(\phi_i,\mu_i,w_i)=\CS(u_i)$ and adjoint states $(p_i,q_i,r_i)$. Then, with a constant $K_8>0$ that depends only on $R$ and the 
data, it holds that
\begin{align}
\label{contdepadj}
&\|(p_1+\tau q_1)-(p_2+\tau q_2)\|_{H^1(0,T;H)\cap C^0([0,T];V)\cap L^2(0,T;\Hdue)} \,+\,\|p_1-p_2\|_{L^2(0,T;H^4(\Omega))}
\nonumber\\
&+\,\|q_1-q_2\|_{L^2(0,T;\Hdue)}\,+\,\|r_1-r_2\|_{H^1(0,T;\Hdue)}\,\le\,K_8\,\|u_1-u_2\|_{L^2(0,T;H)}\,.
\end{align}
\Ecor
\Bdim
We put $p=p_1-p_2,\,\,q=q_1-q_2,\,\,r=r_1-r_2$. Then $(p,q,r)$ is the unique strong solution to the system
\begin{align}
\label{rudi1}
&-\dt(p+\tau q)-\Delta q+f''(\phi_1)q=z_1 &&\mbox{a.e. in }\,Q,\\
\label{rudi2}
&-\Delta p-q=0 &&\mbox{a.e. in }\,Q,\\
\label{rudi3}
&-\gamma \dt r+r-q=0 &&\mbox{a.e. in }\,Q,\\
\label{rudi4}
&\dn p=\dn q=0 &&\mbox{a.e. on }\,\Sigma,\\
\label{rudi5}
&(p+\tau q)(T)=z_2,\quad r(T)=0&&\mbox{a.e. in }\,\Omega, 
\end{align}
where
\begin{equation}
\label{diezs}
z_1=-(f''(\phi_1)-f''(\phi_2))q_2 + b_1(\phi_1-\phi_2)\quad\mbox{and}\quad z_2=b_2(\phi_1(T)-\phi_2(T)).
\end{equation}
Applying the same sequence of estimates that led above in the proof of Theorem~\ref{Teo3} to the derivation of \eqref{adj6} 
(but this time to the continuous system \eqref{rudi1}--\eqref{rudi5}), we readily see that the assertion will be proved
as soon as we can show that
$$
\|z_1\|_{\L2 H}+\|z_2\|_V\,\le\,C\,\|u_1-u_2\|_{\L2 H}\,.
$$
But this is an immediate consequence of the estimate \eqref{contdep2} in Theorem~\ref{Teo2}: indeed, using the continuity of
the embedding $V\subset L^4(\Omega)$, we have that
\begin{align*}
&\|z_1\|^2_{\L2 H}+\|z_2\|^2_V
\\
&\le\,C\int_0^T\|(\phi_1-\phi_2)(s)\|_{L^4(\Omega)}^2\,\|q_2(s)\|^2_{L^4(\Omega)}\,ds
\\
&\qquad 
+ C\|\phi_1-\phi_2\|_{\L2H}^2 \,+\,  C \|(\phi_1-\phi_2)(T)\|_{V}^2
\\
&\le C\,\|\phi_1-\phi_2\|^2_{C^0([0,T];V)}\bigl(1+\|q_2\|_{\L2 V}^2\bigr)\,\le\,C\,\|u_1-u_2\|^2_{\L2 H}.
\end{align*}
\Edim

\subsection{First-order necessary optimality conditions}

In this section, we aim at deriving  associated first-order necessary optimality conditions for local minima of the
optimal control problem {\bf (CP)}. We assume that {\bf (A1)}--{\bf (A5)} are fulfilled and that $\juerg{G}:L^2(0,T;H)\to\erre$ 
is a general \juerg{nonnegative,} convex and continuous functional. We define the reduced cost functionals associated with the functionals  
$J$ and ${\cal J}$  introduced in \eqref{cost} by
\begin{equation}\label{reduced}
\widehat J(u) :=  J(\revis{\CS_1(u)},u), \quad \widehat {\cal J}(u)={\cal J}(\revis{\CS_1(u)},u)\,.
\end{equation}
Since $\CS \revis{{}=(\CS_1, \CS_2,\CS_3)}$ is twice continuously Fr\'echet differentiable  from $\CU$ into the space $C^0([0,T];H)^3$ (which contains $\CZ$), 
it follows from the chain rule that the smooth part $\,\widehat J\,$ of the reduced objective functional is  a twice continuously 
 Fr\'echet differentiable mapping from $\CU$ into $\erre$, where, for every $\us\in\CU$ and every $h\in\CU$, it holds 
with $(\phis,\mus,\ws)=\CS(u^*)$ that
\begin{align}
\label{DJ}
\widehat J'(\us)[h]\,&=\,b_1\iint_Q \xi(\phi^*-\phi_Q)\,+\,b_2\int_\Omega\xi(T)(\phis(T)-\phi_\Omega) 
\,+\, b_3 \iint_Q \us h \,,
\end{align}
where $\,(\xi,\eta,v)=\CS'(\us)[h]\in\CZ\,$ is the unique solution to the linearized system \Linear\
associated with $h$.
\Brem
\label{Rem4.5}
Observe that the \rhs\ of \eqref{DJ} is meaningful also for arguments $h\in L^2(0,T;H)$, where in this case $(\xi,\eta,v)$
is still $\CS'(\us)[h]$, but with the extension of the operator $\CS'(\us)$ to $L^2(0,T;H)$ introduced in Remark~\ref{Spextended}.
 Hence, by means of the
identity \eqref{DJ} we can extend the operator $\widehat J'(\us)\in\CU^*$ to $L^2(0,T;H)$. The extended operator, which we again denote by
$\widehat J'(\us)$, then becomes an element of $L^2(0,T;H)^*$. In this way, expressions of the form $\widehat J'(\us)[h]$
have a proper meaning also for 
$h\in L^2(0,T;H)$.  
\Erem
In the following, we assume that $\us\in\Uad$ is a locally optimal control for ${\bf (CP)}$ in the sense
of \,$\CU$, that is, there is some $\varepsilon>0$ such that
\begin{equation}
\label{lomin}
\widehat {\cal J}(u)\,\ge\,\widehat {\cal J}(\us)\quad\mbox{for all $u\in\Uad$ satisfying }\,\|u-\us\|_{\CU}\,\le\,\varepsilon.
\end{equation}
Notice that any locally optimal control in the sense of $L^p(Q)$ \revis{for some} $1 \le p < \infty$
is also locally optimal in the sense of $\CU$, since the topology of $\CU$ is the finest among these spaces. 
Therefore, a result proved for locally optimal controls in the sense of  $\CU$ is also valid for locally optimal controls 
in the sense of  $L^p(Q)$ for \revis{any} $1\le p<\infty$. It is also true for globally optimal controls.

A standard argument \revis{of nondifferentiable convex optimization theory
(for full details see}, e.g., \cite{SpTr1, SpTr2}) then shows 
that there is 
some $\lambda^* \in \partial \juerg{G}(\us)\subset L^2(0,T;H)$ such that 
\begin{equation} \label{varineq1}
{\widehat J}'(\us)[u - \us] +  \kappa\iint_Q 
\lambda^* (u- \us)\,\ge 0 \quad \forall \,
u \in \Uad.
\end{equation}
As usual, we simplify the expression ${\widehat J}'(\us)[u-\us]$ in \eqref{varineq1}
by means of the adjoint state variables defined in \Adjoint. A standard calculation using the linearized system \Linear\
then leads to the following result.
\Bthm
\,{\rm (Necessary optimality condition)}  \label{Thm4.5}Suppose that {\bf (A1)}--{\bf (A5)}, \revis{{}$\gamma \in W^{2,\infty}(\Omega)$ and{}} $w_0\in L^\infty (\Omega)$
are fulfilled and that \juerg{$G:L^2(0,T;H)\to\erre$ is nonnegative,} convex and continuous. Moreover, 
 let $\us \in \Uad$ be a locally optimal control of  {\bf (CP)} in the sense of $\,\CU\,$ 
with associated state $\,(\phis,\mus,\ws)={\cal S}(\us)$
and adjoint state $\,(\ps,\qs,\rs)$.
Then there exists some $\lambda^*   \in \partial \juerg{G}(\us)$ such that,
for all $\, u \in \Uad$, 
\begin{align}
\label{varineq2}
&\iint_Q  \left(r^*+\kappa\lambda^*+b_3 u^*\right)\left(u-\us\right) \, \ge \,0 \,.
\end{align}
\Ethm
We underline again that \eqref{varineq2} is also necessary for all globally optimal controls and all controls 
which are locally optimal in the sense of $L^p(Q)$ with $p \ge 1$.

\subsection{Sparsity of controls}

The convex function \juerg{$G$} in the objective functional accounts for the sparsity of optimal controls, i.e., the possibility
that any locally optimal control may vanish in some subset of the space-time cylinder $Q$. The form of this region 
depends on the particular choice of the functional \juerg{$G$}.
The sparsity properties can be deduced from the variational inequality \eqref{varineq2} and the particular
form of the subdifferential  $\juerg{\partial G}$. In what follows, we restrict ourselves to the case of {\em full sparsity}
which is connected to the $L^1(Q)-$norm functional \juerg{$G$} introduced in \eqref{defg}. 
Its subdifferential is given by (see \cite{ioffe_tikhomirov1979})
\begin{equation}\label{dg}
\partial G(u) = \Biggl\{\lambda \in L^2(Q):\,
\lambda(x,t) \in \left\{
\begin{array}{ll}
\{1\} & \mbox{ if } u(x,t) > 0\\
{[-1,1]}& \mbox{ if } u(x,t) = 0\\
\{-1\} & \mbox{ if } u(x,t) < 0\\
\end{array}
\right.
\mbox{ for a.e. } {(x,t) \in Q}
\Biggr\}.
\end{equation}
We then have to use this subdifferential in the variational inequality \eqref{varineq2} from Theorem~\ref{Thm4.5} to 
obtain the following result.
\Bthm {\rm (Full sparsity)} \,\,Suppose that the assumptions {\bf (A1)}--{\bf (A5)}, \revis{{}$\gamma \in W^{2,\infty}(\Omega)$ and{}} $w_0\in L^\infty (\Omega)$ are fulfilled, and assume that 
$\underline u$ and $\overline u$ are constants such that $\underline u <0 <\overline u$. Let $\us\in \Uad$ be a 
locally optimal control in the sense 
of \,$\CU$\, for the problem {\bf (CP)}
with the functional $\,\juerg{G}\,$ defined in \eqref{defg}, and with associated state $(\phis,\mus,\ws)=\CS(\us)$ solving \State\ and 
adjoint state $(p^*,q^*,r^*)$ solving \Adjoint. Then there exists a function $\lambda^*\in \juerg{\partial G}(\us)$ 
satisfying \eqref{varineq2}, and it holds
\begin{eqnarray}
\us(x,t) = 0 \quad &\Longleftrightarrow& \quad |r^*(x,t)| \le \kappa, \quad\mbox{for a.e. }\,(x,t)\in Q. 
\label{fullsparse}  
\end{eqnarray}
Moreover, if\, $r^*$ and $\lambda^*$ are given, then
$\us$  is obtained from the projection formula
\begin{eqnarray*}
\us(x,t)& =& \max\left\{\underline u, \min\left\{\overline u, -b_3^{-1} \left(r^*+ \kappa \,\lambda^*\right)(x,t)\right\}\right\}
\,\mbox{ for a.e. $(x,t)\in Q$}.
\end{eqnarray*}
\Ethm
\Bdim
The projection formula is a direct consequence of the variational inequality \eqref{varineq2}.
It remains to show the validity of \eqref{fullsparse}.  
We use the projection formula and the fact that $\underline u<0<\overline u$. For a.e. $(x,t)\in Q$, we have: if $\,\us(x,t)=0$, then
$\,\,-b_3^{-1}(r^*(x,t)+\kappa\lambda^*(x,t))=0$, where $\lambda^*(x,t)\in [-1,1]$. Consequently,
$\,|r^*(x,t)|=\kappa|\lambda^*(x,t)|\le\kappa$.

Now let us assume that $|r^*(x,t)|\le \kappa$. If $\us(x,t)>0$, then $\lambda^*(x,t)=1$ and, by the projection formula,
$\,\,-b_3^{-1}(r^*(x,t)+\kappa)\ge u^*(x,t)>0,$ which implies that $\,r^*(x,t)+\kappa<0\,$ and thus $\,|r^*(x,t)|=-r^*(x,t)>\kappa$,
a contradiction. By analogous reasoning, we can show that also the assumption $\,u^*(x,t)<0$ \,leads to a contradiction. 
We thus must have $u^*(x,t)=0$. 
This ends the proof.
\Edim

We conclude this subsection by showing that all locally optimal controls in the sense of $\CU$  are identically zero
for sufficiently large sparsity parameters. Indeed,
the global estimate \eqref{ssbound1} for the solutions to the state system is valid for all controls $u\in\Uad$, and this 
is also true for the global estimate \eqref{adj6}. Hence, there is some $C^*>0$ such that
\begin{equation*}
\|r^*\|_{H^1(0,T;\Hdue)}\,\le\,C^* \quad\forall\,\us\in\Uad, 
\end{equation*}
and, in view of the continuity of the embedding $H^1(0,T;\Hdue)\subset C^0(\overline Q)$, also
$$\|r^*\|_{C^0(\overline Q)}\,\le\,\kappa^* \quad\forall\,\us\in\Uad,
$$
for a sufficiently large $\kappa^*>0$, which proves our claim.

\subsection{Second-order sufficient optimality conditions} 
We conclude this paper with the derivation of second-order sufficient optimality conditions. 
We provide conditions that ensure local optimality of functions \,$\us\,$ obeying the first-order necessary optimality 
conditions of Theorem~\ref{Thm4.5}. Second-order sufficient optimality conditions are based on a condition of coercivity that 
is required to hold for the smooth part  $\,J\,$ 
of $\,{\cal J}\,$ in a certain critical cone. The nonsmooth part \juerg{$\,G\,$} contributes to sufficiency by its convexity. In the following,
we generally assume that the conditions {\bf (A1)}--{\bf (A5) are} fulfilled.
Our analysis will follow closely the lines of \cite{casas_ryll_troeltzsch2015}, where a second-order analysis was 
performed for sparse control of the FitzHugh--Nagumo system. In particular, we adapt the proof of 
\cite[Thm.~3.4]{casas_ryll_troeltzsch2015} to our setting of a viscous Cahn--Hilliard system.

To this end, we fix  a control $\,\us \,$ that satisfies the first-order necessary optimality conditions,
and we set $\,(\phis,\mus,\ws)=\CS(\us)$. Then the cone
\[
C(\us) = \{ v \in L^2(0,T;H) \,\text{ satisfying the sign conditions \eqref{sign} a.e. in $Q$}\},
\]
where
\begin{equation} \label{sign}
v(x,t) \left\{
\begin{array}{l}
\ge 0 \,\, \text{ if }\,\, u^*(x,t) = \underline u\\ 
\le 0 \,\, \text{ if }\,\, u^*(x,t) = \overline u 
\end{array}
\right. \,,
\end{equation}
is called the {\em cone of feasible directions}, which is a convex and closed subset of $L^2(0,T;H)$.
 We also need the directional derivative of \juerg{$G$} at $u\in L^2(0,T;H)$ in the direction 
$v\in L^2(0,T;H)$, which is given by
\begin{equation}
\label{g'}
\juerg{G'(u,v) = \lim_{t \searrow 0} \frac{1}{t}(G(u+t v)-G(u))}\,.
\end{equation}
Following the definition of the critical cone in \cite[Sect.~3.1]{casas_ryll_troeltzsch2015}, we define
\begin{equation}
\label{critcone}
C_{\us} = \{v \in C(\us): \widehat{J}'(\us)[v] + \kappa \juerg{G}'(\us,v) = 0\}\,,
\end{equation}
which is also a closed and convex subset of $L^2(0,T;H)$. According to \cite[Sect.~3.1]{casas_ryll_troeltzsch2015}, 
it consists of all $v\in C(\us)$ satisfying  
\begin{equation} \label{pointwise}
v(x,t) \left\{
\begin{array}{l}
= 0 \,\,\mbox{ if } \,\,|r^*(x,t) + b_3 u^*(x,t)| \not= \kappa\\
\ge 0 \,\,\mbox{ if }\,\, u^*(x,t) = \underline u\,\, \mbox{ or } \,\,(r^*(x,t) = -\kappa \,\,\mbox{ and }\,\, u^*(x,t) = 0)\\
\le 0 \,\,\mbox{ if }\,\, u^*(x,t) = \overline u\,\, \mbox{ or }\,\, (r^*(x,t) = \kappa \,\,\mbox{ and } \,\,u^*(x,t) = 0)
\end{array}
\right.\,.
\end{equation}

At this point, we derive an explicit expression for $\,\widehat J''(u)[h,k]\,$ 
for arbitrary $\,u,
h,k\in\CU$. In the following, we argue similarly as in 
\cite[Sect.~5.7]{Fredibuch}. At first, we readily infer that, 
for every $((\phi,\mu,w), u)\in (C^0([0,T];H))^3\times \CU$ 
and ${\bf y}=(y_1,y_2,y_3), {\bf z}=(z_1,z_2,z_3)$ such that
$ ({\bf y}, u_1),({\bf z}, u_2)\in (C^0(0,T;H))^3\times \CU$, it follows for the quadratic functional $J$ that
\begin{align}
	\label{D2:0}
	J''((\phi,\mu,w), u)[({\bf y},u_1),({\bf z},u_2)]
	= b_1 \iint_Q  y_1 z_1 \,+\,b_2 \iO y_1(T) z_1(T)
	\,+\,b_3 \iint_Q u_1\,u_2.
\end{align}
For the second-order derivative of the reduced cost functional $\widehat J$
at a fixed control $\us$ we then find with $\,(\phis,\mus,\ws)=\CS(\us)$\, that

\begin{align}
	\nonumber \widehat J''(\us)[h,k] & = D_{(\phi,\mu,w)}J((\phis,\mus,\ws), \us)[(\psi,\nu,z)]\\
	\label{D2:1}
	&\quad + \,J'' ((\phis,\mus,\ws), \us)[((\xi^h,\eta^h,v^h),h),((\xi^k,\eta^k,v^k),k)],
\end{align}
where $(\xi^h,\eta^h,v^h)$, $(\xi^k,\eta^k,v^k)$, and $(\psi,\nu,z)$ stand for 
the unique corresponding solutions to the linearized 
system associated with $h$ and $k$, and to the bilinearized system, respectively.
From the definition of the cost functional \eqref{cost} we readily infer that
\begin{align}
	\label{D2:2} 
	D_{(\phi,\mu,w)} J((\phis,\mus,\ws),\us)[(\psi,\nu,z)] 
	= b_1   \iint_Q (\phis - \phi_Q) \psi
	\,+\, b_2\iO (\phis(T) - \phi_{\Omega})\psi(T).
\end{align}

We now claim that, with the associated adjoint state $\,(p^*,q^*,r^*)$,
\begin{align}
&b_1 \iint_Q (\phis-\phi_Q)\psi 
 	\,+ b_2\iO (\phis(T) - \phi_{\Omega})\psi(T)\,=\,- \iint_Q f^{(3)}(\phis)\xi^h\xi^k q^* \,.
\end{align}
To prove this claim, we multiply \eqref{bilin1} by $\,p^*$, \eqref{bilin2} by $\,q^*$, \eqref{bilin3} by $\,r^*$, 
add the resulting equalities, and integrate over $\,Q$ and by parts, to obtain that
\begin{align*}
	0\,& = \,\iO p^*(T)\psi(T)-\iint_Q \dt p^*\psi\,-\iint_Q \nu\,\Delta p^*\,+\iO \tau q^*(T)\psi(T)
	\,-\iint_Q \tau \dt q^*\psi\nonumber\\
	&\quad -\iint_Q \psi\Delta q^* \,-\iint q^*(\nu+z)\,+\iint f''(\phis)\psi q^*\, + \iint_Q f^{(3)}(\phis)\xih\xik q^*\nonumber\\
	& \quad
	+ \iO \gamma r^*(T) z(T)\,-\iint_Q \gamma \dt r^*z \,+\iint_Q r^*z\nonumber\\ 
	&=\, \iO b_2(\phis(T)-\phi_\Omega)\psi(T)\,+\iint_Q \psi\Big[-\dt(p^*+\tau q^*)-\Delta q^*+f''(\phis)q^*\Big]\,\nonumber\\
	&\quad + \iint_Q \nu\Big[-\Delta p^*-q^*\Big] \,+\iint_Q z\Big[-\gamma\dt r^* +r^*-\juerg{q^*}\Big]\nonumber\\
	&=\,b_1\iint_Q(\phis-\phi_Q)\psi\,+\,b_2\iO(\phis(T)-\phi_\Omega)\psi(T)
	\,+\iint_Q f^{(3)}(\phis)\xih\xik q^*\,,
	\end{align*} 
	whence the claim follows, since $\,(p^*,q^*,r^*)\,$ solves the adjoint system \eqref{adj1}--\eqref{adj5}.
From this characterization, along with \eqref{D2:1} and \eqref{D2:2}, we conclude that
\begin{align}	
& \widehat J''(\us)[h,k] \,=\,\iint_Q \bigl(b_1-f^{(3)}(\phis)q^*\bigr)\xih\,\xik \,+\,b_2\iO \xih(T)  \xik(T)
\,+\,b_3\iint_Q h \,k\,.
\label{walter1}
\end{align}

Observe that the expression on the right-hand side of \eqref{walter1} is meaningful also for increments $\,h,k\in L^2(Q)$. 
Indeed, in this case the expressions $\,(\xih,\eta^h,v^h)=\CS'(\us)[h]$, $\,(\xi^k,\eta^k,v^k)=\CS'(\us)[k]$,
and $\,(\psi,\nu,z)=\CS''(\us)[h,k]\,$ have an interpretation in the sense of the extended operators $\,\CS'(\us)\,$ and
$\,\CS''(\us)\,$ introduced in Remark~\ref{Spextended}. Therefore, the operator $\,\widehat J''(\us)\,$ can be extended by
the identity \eqref{walter1} to the space $L^2(Q)\times L^2(Q)$. This extension,
which will still be denoted by $\,\widehat J''(\us)$, will be frequently used in the following.
We now show that it is continuous. Indeed, we claim 
that for all $h,k\in L^2(Q)$ it holds 
\begin{align}
\label{walter2}
\left|\widehat J''(\us)[h,k]\right|\,\le\,\widehat C\,\|h\|_{L^2(Q)}\,\|k\|_{L^2(Q)}\,,
\end{align}
where the constant $\widehat C>0$ is independent of the choice of $\,\us \in {\cal U}_R$. Obviously, only the first integral on the 
right-hand side
of \eqref{walter1} needs some treatment. We have,
by virtue of H\"older's inequality, the continuity of the embedding $V\subset L^4(\Omega)$, and the global bounds \eqref{ssbound2}, 
\eqref{extension}, and \eqref{adj6}, 
\begin{align*}
&\Big| \iint_Q f^{(3)}(\phis) \xih\xik q^* \Big|\,\le\,C\int_0^T \|\xih(t)\|_{L^4(\Omega)}\,\|\xik(t)\|_{L^4(\Omega)}\,
\|q^*(t)\|_{L^2(\Omega)}\,dt\\[2mm]
&\quad\le \,C\,\|\xih\|_{C^0([0,T];V)}\,\|\xik\|_{C^0([0,T];V)}\,\|q^*\|_{L^2(0,T;H)}\,\le\,
C\,\|h\|_{L^2(Q)}\,\|k\|_{L^2(Q)}\,, 
\end{align*}
as asserted. 

In the following, we will employ the following coercivity condition:

\begin{equation} \label{coerc}
\widehat J''(\us)[v,v] >0 \quad \forall \, v \in C_{\us} \setminus \{0\}\,.
\end{equation}

\noindent Condition \eqref{coerc} is a direct extension of associated conditions that are standard in finite-dimensional 
nonlinear optimization. In the optimal control of partial differential equation, it was first used in \cite{casas_troeltzsch2012}.
We have the following result.
\begin{theorem} \,\,{\rm (Second-order sufficient condition)} \,\,Suppose that {\bf (A1)}--{\bf (A5)} are fulfilled 
along with \revis{{}$\gamma \in W^{2,\infty}(\Omega)$ and{}} $w_0\in L^\infty (\Omega)$. 
Moreover, let $\us \in \Uad$, together with the associated state 
$(\phi^*,\mu^*,w^*)=\CS(\us)$ and the adjoint state $(p^*,q^*,r^*)$, fulfill the first-order necessary optimality conditions 
of Theorem~\ref{Thm4.5}. If, in addition,  $\us$ satisfies the coercivity condition \eqref{coerc}, 
then there exist constants $\,\varepsilon > 0\,$ and $\,\zeta > 0\,$ such that the quadratic growth condition
\begin{equation} \label{growth}
\widehat {\cal J}(u) \ge \widehat {\cal J}(\us) + \zeta \, \|u-\us\|^2_{L^2(Q)} 
\end{equation}
holds for all $\,u \in \Uad\,$ with $ \|u-\us\|_{L^2(Q)}  < \varepsilon$. Consequently,
$\,\us\,$ is a locally optimal control in the sense of $\,L^2(Q)$.
\end{theorem}
\begin{proof} The proof follows that of \cite[Thm.~3.4]{casas_ryll_troeltzsch2015}. We include it for the reader's convenience.
We argue by contradiction, assuming that the claim of the theorem is not true. Then there exists a sequence of controls 
$\{u_j\}\subset \Uad$ such that, for all $j\in\enne$, 
\begin{equation}
\label{contrary}
\|u_j-\us\|_{L^2(Q)} < \frac{1}{j} \quad \mbox{ while } \quad \widehat {\cal J}(u_j)<  \widehat {\cal J}(\us) 
+ \frac{1}{2j} \|u_j-\us\|_{L^2(Q)}^2\,.
\end{equation}
Noting that $u_j \not=\us$ for all $j\in\enne$, we define 
$$
\tau_j := \|u_j-\us\|_{L^2(Q)}
  \quad \mbox{ and } \quad h_j := \frac{1}{\tau_j}(u_j-\us)\,.
$$                                               
Then $\|h_j\|_{L^2(Q)}=1$ and, possibly after selecting a subsequence, we can assume that 
\[
h_j \to h \, \mbox{ weakly in }\,  L^2(Q)
\]
for some $h\in L^2(Q)$. As in \cite{casas_ryll_troeltzsch2015}, the proof is split into three parts. 

(i) $h \in C_{\us}$: Obviously, each $h_j$ obeys the sign conditions \eqref{sign} and thus belongs to $C(\us)$. 
Since $C(\us)$ is convex and closed in $L^2(Q)$, it follows that $h\in C(\us)$. We now claim that 
\begin{equation}
\label{Paulchen}
\widehat J'(\us)[h] + \kappa \juerg{G}'(\us,h) = 0.
\end{equation}
Notice that by Remark~\ref{Rem4.5} the expression 
$\,\widehat{J}'(\us)[h]\,$ has a well-defined meaning. For every $\vartheta \in (0,1)$ and all $h,u \in L^2(Q)$, 
we infer from the convexity of $\,\juerg{G}\,$ that
\begin{align}
\juerg{G}(h)-\juerg{G}(u) &\,\ge\, \frac{\juerg{G}(u + \vartheta (h-u))-\juerg{G}(u)}{\vartheta} \,\ge\, \juerg{G}'(u,h-u)
\,=\, \max_{\lambda \in \partial \juerg{G}(u)}\,
\iint_Q \lambda (h-u),
\label{directionalder}
\end{align}
where the last equality can be checked directly using \eqref{defg}, \eqref{g'} and \eqref{dg}. Then,
in particular, we have that
\begin{align}
&\widehat{J}'(\us)[h] +  \kappa \juerg{G}'(\us,h)\,\ge\, \widehat{J}'(\us)[h] + \iint_Q \kappa \lambda^* h 
= \iint_Q  (r^* + b_3 u^*+\kappa\lambda^*) h \nonumber\\
&\quad= \lim_{j \to \infty} \frac{1}{\tau_j}  \iint_Q (r^* + b_3 u^* + \kappa \lambda^*) (u_j-u^*)\big)
\,\ge\, 0\,,
\end{align}
by the variational inequality \eqref{varineq2}. Next, we prove the converse inequality. By \eqref{contrary}, 
it turns out that
\[\widehat{J}(u_j)-\widehat{J}(\us) + \kappa\left(\juerg{G}(u_j)-\juerg{G}(\us)\right) < \frac{1}{2j} \tau_j^2\,,
\]
whence, owing to the mean value theorem, and since $u_j = \us+  \tau_j h_j$, we get
\[
\tau_j \widehat{J}'(\us + \theta_j \tau_j h_j)[h_j] + \kappa (\juerg{G}(\us + \tau_j h_j)-\juerg{G}(\us) )
<  \frac{1}{2j} \tau_j^2\,, 
\]
with some \juerg{real number $\theta_j \in (0,1)$}. 
Now observe that the mapping $h\mapsto \juerg{G}'(\us,h)$ is positive homogeneous on $L^2(Q)$.  
We therefore obtain from \eqref{directionalder} that 
$$
\kappa(\juerg{G}(\us + \tau_j h_j)-\juerg{G}(\us))\ge \kappa\, \juerg{G}'(\us,\tau_j h_j) =\tau_j\,\kappa\,\juerg{G}'(\us,h_j),
$$
so that, after division by $\,\tau_j$,
\begin{equation}
\label{Yeah}
\widehat{J}'(\us + \theta_j \tau_j h_j)[h_j] + \kappa \juerg{G}'(\us,h_j) < \frac{ \tau_j}{2j}\,.
\end{equation}
At this point, we note that the mapping $h\mapsto \juerg{G}'(\us,h)$ is convex and continuous, and thus weakly sequentially 
semicontinuous, on $L^2(Q)$. Consequently, 
$$
\juerg{G}'(\us,h)\,\le\,\liminf_{j\to\infty} \,\juerg{G}'(\us,h_j).
$$
Besides, the sequence $\{\widetilde u_j\}_{j\in\enne}\subset\Uad$, where $\,\widetilde u_j=\us+\theta_j \tau_j h_j
=\us+\theta_j(u_j-\us)$, converges strongly in $L^2(Q)$ to $\us$. Now let, for $j\in\enne$,
 $(\widetilde \phi_j,\widetilde \mu_j,\widetilde w_j)
=\CS(\widetilde u_j)$, and let $(\widetilde p_j,\widetilde q_j,\widetilde r_j)$ denote the associated adjoint state.
Then, by Corollary~\ref{Cor4.4}, $\widetilde r_j\to r^*$ strongly in $H^1(0,T;\Hdue)$ and thus $\widetilde r_j+b_3\widetilde u_j
\to r^*+b_3 \us$ strongly in $L^2(Q)$, as $j\to\infty$. Consequently, we have that
\begin{align*}
\lim_{j\to\infty} \widehat J'(\widetilde u_j)[h_j]=\lim_{j\to\infty}\iint_Q (\widetilde r_j+b_3\widetilde u_j)h_j
=\iint_Q(r^*+b_3\us)h=\widehat J'(\us)[h],
\end{align*}
and we obtain from \eqref{Yeah} that
\begin{align*}
&\widehat{J}'(\us)[h] + \kappa \juerg{G}'(\us,h) \le \lim_{j\to\infty} \widehat J'(\widetilde u_j)[h_j]+\kappa
\liminf_{j\to\infty} \juerg{G}'(\us,h_j)\\
&=\,\liminf_{j\to\infty}\,\bigl(\widehat J'(\widetilde u_j)[h_j]+ \kappa \juerg{G}'(\us,h_j)\bigr)\,\le\,0,
\end{align*}
which completes the proof of (i).

(ii) $h = 0$: We again invoke \eqref{contrary}, now  performing a second-order Taylor expansion on the left-hand side\juerg{. We
obtain, with some real number $\theta_j\in (0,1)$,}
\begin{align}
&\widehat{J}(\us) + \tau_j  \widehat{J}'(\us)[h_j] + \frac{\tau_j^2}{2} \widehat{J}''(\us + \theta_j \tau_j h_j)[h_j,h_j]
+ \kappa \juerg{G}(\us + \tau_j h_j)\nonumber \\
&\quad<\widehat{J}(\us) + \kappa \juerg{G}(\us) +  \frac{ \tau_j^2}{2j}\,.\nonumber 
\end{align}
We subtract $\,\widehat{J}(\us) + \kappa \juerg{G}(\us)\,$ from both sides and use \eqref{directionalder} once more to find that
\begin{equation} \label{Konrad2}
\tau_j  \left(\widehat{J}'(\us)[h_j]  + \kappa \juerg{G}'(\us,h_j)\right)+ \frac{\tau_j^2}{2} \widehat{J}''(\us + \theta_j \tau_j h_j)
[h_j,h_j]< \frac{ \tau_j^2}{2j}\,.
\end{equation}
From the right-hand side of \eqref{directionalder} and the variational inequality \eqref{varineq2}, it follows that
\[
\widehat{J}'(\us)[h_j]  + \kappa \juerg{G}'(\us,h_j) \ge 0\,,
\]
and thus, by \eqref{Konrad2},
\begin{equation}\label{liminf1} 
\widehat{J}''(\us + \theta_j \tau_j h_j)[h_j,h_j]< \frac{1}{j}\,.
\end{equation}
At this point, we apply the identity \eqref{liminf0} in Lemma~\ref{Lem4.9} below, where we note that $\widetilde u_j=\us+\theta_j\tau_j h_j$
converges to $\us$ strongly in $\L2 H$. We have, using the notation introduced in Lemma~\ref{Lem4.9} and \eqref{liminf1},
\begin{align*}
&\widehat J''(\us)[h,h]\,=\,\lim_{j\to\infty} \Big(\iint_Q\bigl(b_1-f^{(3)}(\widetilde\phi_j)
\widetilde q_j\bigr)\bigl|\widetilde\xi^{h_j}\bigr|^2 \,+\,b_2\iO \bigl|\widetilde\xi^{h_j}(T)\bigr|^2 \Big) 
\,+\,b_3\iint_Q|h|^2\\
&\le\,\liminf_{j\to\infty} \Big(\iint_Q\bigl(b_1-f^{(3)}(\widetilde\phi_j)
\widetilde q_j\bigr)\bigl|\widetilde\xi^{h_j}\bigr|^2 \,+\,b_2\iO \bigl|\widetilde\xi^{h_j}(T)\bigr|^2  
\,+\,b_3\iint_Q|h_j|^2\Big)\\
&=\,\liminf_{j\to\infty} \,\widehat J''(\widetilde u_j)[h_j,h_j]\,\le\,0\,. 
\end{align*}
Since we know that $h \in C_{\us}$, the second-order condition \eqref{coerc} implies that \,$h = 0$. 

(iii) {\em Contradiction:} 
From the previous step we know that $\,h_j\to 0\,$ weakly in $\,L^2(Q)$. Now, \eqref{walter1} yields that
\begin{align}
\label{umba-umba}
&\widehat J_1''(\us)[h_j,h_j] = \iint_Q \bigl(b_1-f^{(3)}(\us)q^*\bigr)|\xi^{h_j}|^2\,+\,b_2\iO|\xi^{h_j}(T)|^2
\,+\,b_3\iint_Q |h_j|^2\,,
\end{align}
where we have set $\,(\xi^{h_j},\eta^{h_j},v^{h_j})=\CS'(\us)[h_j]$, for $\,j\in\enne$. Since $h_j\to 0$ weakly in
$\L2 H$, we find from \eqref{liminf0} in Lemma~\ref{Lem4.9} that the sum of the first two integrals on the right-hand side of 
\eqref{umba-umba} converges to zero. 
On the other hand, $\,\|h_j\|_{L^2(Q)}=1$ for all
$\,j\in\enne$, by construction. The weak sequential semicontinuity of norms then implies that 
\begin{align}
\label{uff-uff}
&\liminf_{j\to\infty} \,\widehat J''(\us)[h_j,h_j] \,\ge \,\liminf_{j\to\infty} \,b_3\iint_Q|h_j|^2\,=\,
b_3\,>0\,.
\end{align}
On the other hand, we may apply Lemma~\ref{Lem4.9} twice, namely to the sequence $\widetilde u_j=\us+\theta_j\tau_j h_j$ and 
to the constant sequence $\widetilde u_j=\us$, to infer that 
$$
\lim_{j\to\infty} (\widehat J''(\us)-\widehat J''(\us+\theta_j\tau_j h_j))
[h_j,h_j]=0.
$$
Therefore, thanks to \eqref{liminf1} it is clear that
\begin{align*}
&\liminf_{j \to \infty} \widehat{J}''(\us)[h_j,h_j] \\
&=\,\liminf_{j\to\infty} \Big((\widehat J''(\us)-\widehat J''( 
\us+\theta_j\tau_j h_j))
[h_j,h_j]+\widehat J''(\us+\theta_j\tau_j h_j)[h_j,h_j]\Big)\\
&\le \,\liminf_{j\to\infty} \Big((\widehat J''(\us)-\widehat J''(\us+\theta_j\tau_j h_j))
[h_j,h_j]+\,1/j\Big) \,=\,0\,,
\end{align*}
which contradicts \eqref{uff-uff}. The assertion of the theorem is thus proved. 
\end{proof}

We conclude the paper with the auxiliary result that was used in the above proof.
\begin{lemma} 
\label{Lem4.9}
Assuume that {\bf (A1)}--{\bf (A5)}, \revis{{}$\gamma \in W^{2,\infty}(\Omega)$ and{}} $w_0\in L^\infty (\Omega)$ are satisfied.
Suppose that $\{\widetilde u_j\}\subset\Uad$ converges strongly in $L^2(0,T;H)$ to $\us\in\Uad$, and that $\{h_j\}\subset L^2(Q)$
converges weakly in $L^2(Q)$ to $h$. In addition, let $(\widetilde\phi_j,\widetilde\mu_j,\widetilde w_j)=
\CS(\widetilde u_j)$, and let $\,(\widetilde p_j,\widetilde q_j,\widetilde r_j)$ be the associated adjoint state. Moreover, let, 
for arbitrary $h\in \L2 H$, $(\xih,\eta^h,v^h)=\CS'(\us)[h]$, as well as $(\widetilde \xi^{h_j},\widetilde\eta^{h_j},
\widetilde v^{h_j})=\CS'(\widetilde u_j)[h_j]$.
 Then
\begin{align} \label{liminf0}
&\lim_{j \to \infty} \Big(\iint_Q\bigl(b_1-f^{(3)}(\widetilde\phi_j)\widetilde q_j\bigr)\bigl|\widetilde\xi^{h_j}\bigr|^2
\,+\,b_2\iO \bigl|\widetilde\xi^{h_j}(T)\bigr|^2 \Big) \nonumber\\
&\quad =\iint_Q\bigl(b_1-f^{(3)}(\phis)q^*\bigr)\bigl|\xi^h\bigr|^2\,+\,b_2\iO\bigl|\xi^h(T)\bigr|^2\,.
\end{align} 
\end{lemma}
\begin{proof}
At first, notice that
$$
(\widetilde\xi^{h_j}, \widetilde\eta^{h_j}, \widetilde v^{h_j})-(\xi^h,\eta^h,v^h) = \left(\CS'(\widetilde u_j)
-\CS'(\us)\right)[ h_j]\, +\,\CS'(\us)[h_j-h]\,.
$$
By virtue of \eqref{lip1} (recall Remark~\ref{Spextended} in this regard) and the boundedness of $\{h_j\}$ in $L^2(0,T;H)$, the 
first summand on the right converges strongly to zero in $\,\CX$.
The second converges to zero weakly star in (cf.~\eqref{defX})
$$ \bigl(H^1(0,T;H)\cap \L\infty V\cap L^2(0,T;W)\bigr) \times L^2(0,T;W) \times H^1(0,T;H). $$
Thanks to the compact embeddings $V \subset L^p(\Omega)$ for $1\le p<6$, the compactness result stated in 
\cite[Sect.~8, Cor.~4]{Simon}) ensures that
\begin{equation}
\label{Hugo3}
\widetilde\xi^{h_j}\to \xi^h \quad\mbox{strongly in }\,C^0([0,T];L^5(\Omega))\,.
\end{equation}
In particular, we have that
\begin{align} \label{Hugo4}
&\lim_{j\to \infty} \Big(b_1\iint_Q \bigl|\widetilde\xi^{h_j}\bigr|^2\,+\,b_2\iO\bigl|\widetilde\xi^{h_j}(T)\bigr|^2 \Big)
= b_1 \iint_Q \bigl|\xi^h\bigr|^2 \,+\,b_2\iO\bigl|\xi^h(T)\bigr|^2\,.
\end{align}
Moreover, we obtain from \eqref{contdep2} that $\,\|\widetilde \phi_j-\phis\|_{C^0([0,T];V)}\to 0$, so 
that we can conclude from the
global estimate \eqref{ssbound2} and \eqref{contdepadj} that, as $j\to \infty$,
\begin{align}
\label{Hugo5}
&\|f^{(3)}(\widetilde\phi_j)-f^{(3)}(\phis)\|_{C^0([0,T];L^5(\Omega))} \,\to\,0,\\
\label{Hugo6}
&\|\widetilde q_j-q^*\|_{L^2(0,T;\Hdue)}\,\to\,0.
\end{align}
Combining this with \eqref{Hugo3}, we readily verify that
\begin{equation}
\lim_{j\to\infty}\iint_Q f^{(3)}(\widetilde\phi_j)\widetilde q_j\bigl|\xi^{h_j}\bigr|^2\,=\,\iint_Q
f^{(3)}(\phis)q^*\bigl|\xi^h\bigr|^2,
\end{equation} 
which concludes the proof.
\end{proof}


\smallskip

\revis{%
\section*{Acknowledgments}
This research has been supported by the MIUR-PRIN Grant
2020F3NCPX ``Mathematics for industry 4.0 (Math4I4)''. 
In addition, PC gratefully acknowledges his affiliation 
to the GNAMPA (Gruppo Nazionale per l'Analisi Matematica, 
la Probabilit\`a e le loro Applicazioni) of INdAM (Isti\-tuto 
Nazionale di Alta Matematica).}

\bigskip


\End{document}

\section*{Statements and Declarations}
\subsection*{Funding}
The authors declare that no funds \revis{or grants, other than those mentioned in the Acknowledgements,}
were received during the preparation of this manuscript. 
\subsection*{Competing interests}
The authors have no relevant financial or non-financial interests to disclose.
\subsection*{Author Contributions}
All authors contributed to the paper, in particular to  material preparation and analysis. They all read and approved the final manuscript.
